\documentclass{article}
\usepackage{amssymb}
\usepackage{amsmath}

\usepackage[all]{xypic}

\usepackage{geometry}
\input{xy}
\xyoption{all}

\newcommand{\A}{\scriptscriptstyle{A}}
\newcommand{\B}{\scriptscriptstyle{B}}
\newcommand{\D}{\scriptscriptstyle{D}}
\newcommand{\C}{\scriptscriptstyle{\mathcal{C}}}
\newcommand{\smD}{\scriptscriptstyle{\mathcal{D}}}
\newcommand{\E}{\scriptscriptstyle{E}}
\newcommand{\F}{\scriptscriptstyle{F}}
\newcommand{\G}{\scriptscriptstyle{G}}
\newcommand{\sH}{\scriptscriptstyle{H}}

\newcommand{\J}{\scriptscriptstyle{J}}
\newcommand{\K}{\scriptscriptstyle{K}}
\newcommand{\sL}{\scriptscriptstyle{L}}
\newcommand{\sfL}{\scriptscriptstyle{\mathfrak{L}}}
\newcommand{\sP}{\scriptscriptstyle{P}}
\newcommand{\Q}{\scriptscriptstyle{Q}}
\newcommand{\M}{\scriptscriptstyle{M}}
\newcommand{\sfM}{\scriptscriptstyle{\mathfrak{M}}}

\newcommand{\R}{\scriptscriptstyle{R}}
\newcommand{\sfR}{\scriptscriptstyle{\mathfrak{R}}}
\newcommand{\sS}{\scriptscriptstyle{S}}
\newcommand{\T}{\scriptscriptstyle{T}}
\newcommand{\U}{\scriptscriptstyle{U}}
\newcommand{\V}{\scriptscriptstyle{V}}
\newcommand{\W}{\scriptscriptstyle{W}}

\newcommand{\Cat}{\mathbf{Cat}}
\newcommand{\Mnd}{\mathbf{Mnd}}
\newcommand{\Adj}{\mathbf{Adj}}

\begin{document}

\title{Hopf Parametric Adjoint Objects through a 2-adjunction of the type Adj-Mnd}
\author{\\Adrian Vazquez-Marquez\footnote{Contact: avazquez@uiwbajio.mx} \\ \\
  Universidad Incarnate Word, Campus Baj\'io}
\date{}

\maketitle

\begin{center}
\section*{Abstract}
\end{center}

\emph{In this article Hopf parametric adjunctions are defined and analysed within the context of the 2-adjunction of the type $\Adj$-$\Mnd$. In order to do so, the definition of adjoint objects in the 2-category of adjunctions and in the 2-category of monads for $Cat$ are revised and characterized. This article finalises with the application of the obtained results on current categorical characterization of Hopf Monads.}\\

$\phantom{br}$

\begin{flushright}
\emph{In memory of Lecter and Cosmo}
\end{flushright}

\section*{Introduction}

In 2002, I. Moerdijk \cite{moie_motc} characterized the liftings of a monoidal structure to the category of Eilenberg-Moore algebras, for a related initial monad. This characterization lead to the definition of a opmonoidal monad. In 2011, A. Brugui\`eres, S. Lack and A. Virelizier \cite{bral_homm} characterized the liftings of a closed monoidal structure through the concept of a Hopf monad. These two examples will be analysed in the context of higher category theory.
\\

This article belong to a series where 2-adjunctions of the type $\Adj$-$\Mnd$ are applied to \emph{classical} monad theory. In this installment, the author analyse adjoint objects and parametric adjunctions within this context.\\

On the last subject, that of parametric adjunctions, this article is mainly based upon the ideas laid out in the seminal article of A. Brugui\`eres \emph{et. al.} \cite{bral_homm}. In order to apply the 2-adjunction $\Adj$-$\Mnd$, the ideas are developed into a 2-categorical framework, cf. \cite{lojt_apke}. It is in this framework that the Hopf monad concept, for a monoidal closed structure, is extended to the concept of Hopf 1-cells and adjoint parametric objects on certain 2-categories.\\

Without any further ado, the structure of the article is given.\\

In chapter 1, the 2-categorical structure needed for the rest of the article is given, \emph{namely} the construction of the 2-adjunction $\Phi_{\E}^{\sfM}\dashv\Psi_{\E}^{\sfM}$.\\

In chapter 2, adjoints objects in the 2-categories $\Adj_{\R}(Cat)$ and $\Mnd(Cat)$ are revised and characterized. The characterization of such objects is done based on \cite{bral_homm} which is suitable for the 2-categorical context of the article.\\

In chapter 3, the concept of (left) Hopf 1-cells is defined within the 2-category $\Adj_{\R}(Cat)$ and it is used in order to construct the Hopf parametric adjoint objects in that 2-category. The concept of antipode is reconstructed there.\\

In chapter 4, the concept of Hopf 1-cells is provided within the 2-category $\Mnd(Cat)$ and it is used in order to construct the corresponding Hopf parametric adjoint objects for that 2-category.\\

In chapter 5, the condition for being Hopf 1-cell and the related struture, that of a parametric adjoint object, is analysed through the 2-adjunction. The condition for a 1-cell to be Hopf 1-cell is preserved and then Hopf 1-cells, in each 2-category, are compared using the 2-adjunction. At the end, using the isomorphism of categories, induced by the 2-adjunction, a bijection of Hopf parametric adjoint objects is found.\\

In chapter 6, remaining concepts and statements are done in order to get the main theorem of the article which gives a bijection between Hopf parametric liftings, mimicking those liftings for the closure of a monoidal structure \cite{bral_homm}, and certain parametric adjoint objects is given. This chapter finalises with the corresponding application to Hopf monads in a monoidal category which was the main inspiration for this extension to a 2-categorical context.\\

A list of the notations and conventions taken in this article is given as follows. Consider an adjunction $\mathfrak{L}\dashv\mathfrak{R}$, its unit and counit are denoted as $\eta^{\sfR\sfL}$ and $\varepsilon^{\sfL\sfR}$, respectively. This notation might be complicated but refrain one from running out of, and \emph{a posteriori} very needed, the finite set of greek letters. Nevertheless, the notation will be simplyfied whenever possible. For example, if the adjunction comes from a free-forgetful case, \emph{i.e.} $F^{\sS}\dashv U^{\sS}$, the unit can be written as $\eta^{\U\!\F\!\sS}$ or when a parametric adjunction is involved, on $P$, $F_{\sP}\dashv G_{\sP}$ its unit can be written as $\eta^{\G\!\F\!\sP}$. The direction of the adjunction $\mathfrak{L}\dashv\mathfrak{R}$ will be taken as the direction of its left adjoint functor $\mathfrak{L}$, therefore the domain category of the adjunction is the domain category of $\mathfrak{L}$. The triangular identity given by $\varepsilon^{\sfL\sfR}\mathfrak{L}\circ \mathfrak{L}\hspace{0.5pt}\eta^{\sfR\sfL} = 1_{\sfL}$ will be refered to as the triangular identity associated to $\mathfrak{L}$.\\

For the 2-category $\Cat$, of small categories and functors, the notation $\mathfrak{C}$ will be used instead.\\

The notation $1_{\ast}$  will be used for cases like $1_{\mathcal{P}^{\E}}$ whenever the context allows it. Also, in the cartesian monoidal structure for $Cat$, whenever possible $L\times\mathcal{P}$ will be understood as $L\times 1_{\mathcal{P}}$, for example.\\

In the proof of the communativity for a diagram, arguments based on the naturality property of a certain transformation will ommited whenever possible, in order to spare for the numerous details. The pasting composition of 2-cells will be denoted as $\cdot_{p}$.

\section{The 2-category context}

The 2-category context needed for this article is the following 2-adjunction

\begin{equation*}
\Phi_{\E}^{\sfM}\dashv\Psi_{\E}^{\sfM}:\Adj_{\R}(\Cat) \longrightarrow \Mnd(\Cat)
\end{equation*}

$\phantom{br}$

\noindent the subindex $E$ refers to the Eilenberg-Moore objects, since $\Cat$ admits the construction of algebras, and the superindex $\mathfrak{M}$ refers to the monad case \cite{lojt_apke}. Whenever possible, one or both indexes will be dropped. 

\subsection{The 2-category $\Adj_{\R}(\mathfrak{C})$}


The n-cells for the domain 2-category, $\Adj_{\R}(\mathfrak{C})$, are the following:

\begin{itemize}

\item [i)] The 0-cells are adjunctions $L\dashv R:\mathcal{C}\longrightarrow\mathcal{X}$.

\item [ii)] The 1-cells are of the form $(J,V, \:\lambda^{\J\V}\!, \:\rho^{\J\V}):L\dashv R \longrightarrow \overline{L}\dashv\overline{R}$ and depicted as

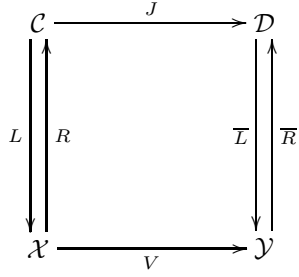
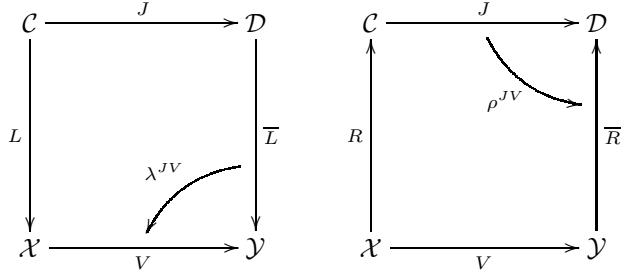
\begin{equation*}
\begin{array}{ccc}
\xy<1cm,1cm>
\POS (0,0) *+{\mathcal{C}} = "a11",
\POS (30,0) *+{\mathcal{D}} = "a12",
\POS (0,-30) *+{\mathcal{X}} = "a21",
\POS (30,-30) *+{\mathcal{Y}} = "a22",
\POS (28,-19) *+{} = "x1",
\POS (15,-29) *+{} = "x2",
\POS "a11" \ar^{J} "a12",
\POS "a12" \ar@<-3pt>_{\overline{L}} "a22",
\POS "a22" \ar@<-3pt>_{\overline{R}} "a12",
\POS "a11" \ar@<-3pt>_{L} "a21",
\POS "a21" \ar@<-3pt>_{R} "a11",
\POS "a21" \ar_{V} "a22",
\endxy 
\end{array}
\end{equation*}

\noindent where

\begin{equation*}
\begin{array}{ccc}
\xy<1cm,1cm>
\POS (0,0) *+{\mathcal{C}} = "a11",
\POS (30,0) *+{\mathcal{D}} = "a12",
\POS (0,-30) *+{\mathcal{X}} = "a21",
\POS (30,-30) *+{\mathcal{Y}} = "a22",
\POS (29,-19) *+{} = "x1",
\POS (15,-29) *+{} = "x2",
\POS "a11" \ar^{J} "a12",
\POS "a12" \ar^{\overline{L}} "a22",
\POS "a11" \ar_{L} "a21",
\POS "a21" \ar_{V} "a22",
\POS "x1"  \ar@/_0.6pc/_{\lambda^{\J\V}} "x2",    
\endxy & &
\xy<1cm,1cm>
\POS (0,0) *+{\mathcal{C}} = "a11",
\POS (30,0) *+{\mathcal{D}} = "a12",
\POS (0,-30) *+{\mathcal{X}} = "a21",
\POS (30,-30) *+{\mathcal{Y}} = "a22",
\POS (15,-1) *+{} = "x1",
\POS (29,-11) *+{} = "x2",
\POS "a21" \ar^{R} "a11",
\POS "a11" \ar^{J} "a12",
\POS "a21" \ar_{V} "a22",
\POS "a22" \ar_{\overline{R}} "a12",
\POS "x1"  \ar@/_0.6pc/_{\rho^{\J\V}} "x2",  
\endxy
\end{array}
\end{equation*}

\noindent are mates and such that $\rho^{\J\V}$ is an isomorphism. The inverse of $\rho^{\J\V}$ will be denoted as $\delta^{\J\V}$ or $\varrho^{\J\V}$. Because of the previous, the notation can be shorten to $(J, V, \lambda^{\J\V})$ or even to $(J,V):L\dashv R \longrightarrow \overline{L}\dashv\overline{R}$, whenever the \emph{left} mate is understood or unimportant. Since the \emph{right} mate is an isomorphism, the 2-category will be denoted as $\Adj_{\R}(\mathfrak{C})$.\\

Note: In general, the mate of a natural transformation $\vartheta:\overline{L}F\longrightarrow GL$ might be denoted as

\begin{equation*}
{}_{\R}\mathbf{m}_{\overline{\R}}\hspace{0.5pt}(\vartheta) = \varepsilon^{\sL\!\R}\cdot_{p}\vartheta \cdot_{p} \eta^{\overline{\R\sL}} = \overline{R}G\varepsilon^{\sL\R}\circ \overline{R}\vartheta R\circ \eta^{\overline{\R\sL}}FR:FR\longrightarrow \overline{R}G.
\end{equation*}

\item [iii)] The 2-cells are comprissed of a pair of natural transformations $(\alpha,\beta)$ where $\alpha:J\longrightarrow J'$ and $\beta:V\longrightarrow V'$ such that one of the following equivalent requirements is fulfilled

\begin{itemize}

\item [$\dagger$)] $\beta L\circ\lambda^{\J\V} = \lambda^{\J'\V'}\circ\overline{L}\alpha $

\item [$\ddagger$)] $\overline{R}\beta\circ \rho^{\J\V} = \rho^{\J'\V'}\circ\alpha R$
  
\end{itemize}

\end{itemize}


\subsection{The 2-category $\Mnd(\mathfrak{C})$}


The $n$-cells for the 2-category $\Mnd(\mathfrak{C})$ are described as follows.

\begin{itemize}

\item [i)] The 0-cells are monads $(\mathcal{C}, S, \mu^{\sS}, \eta^{\sS})$, whose short notation is $(\mathcal{C}, S)$.

\item [ii)] The 1-cells are pairs of the form $(B, \psi^{\B}):(\mathcal{C}, S)\longrightarrow (\mathcal{D}, T)$

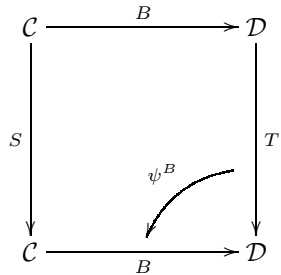
\begin{equation*}
\xy<1cm,1cm>
\POS (0,0) *+{\mathcal{C}} = "a11",
\POS (30,0) *+{\mathcal{D}} = "a12",
\POS (0,-30) *+{\mathcal{C}} = "a21",
\POS (30,-30) *+{\mathcal{D}} = "a22",
\POS (28,-19) *+{} = "x1",
\POS (15,-29) *+{} = "x2",
\POS "a11" \ar^{B} "a12",
\POS "a12" \ar^{T} "a22",
\POS "a11" \ar_{S} "a21",
\POS "a21" \ar_{B} "a22",
\POS "x1" \ar@/_0.6pc/_{\psi^{\B}} "x2",
\endxy
\end{equation*}

\noindent where this natural transformation fulfills the following equations

\begin{eqnarray*}
\psi^{\B}\circ\mu^{\T}B &=& B\mu^{\sS}\circ \psi^{\B}S \circ T\psi^{\B}\\
\psi^{\B}\circ\eta^{\T}B &=& B\eta^{\sS}
\end{eqnarray*}

\noindent these equations might be refered to as the \emph{compatibility}, with the product and the unit of the monads, \emph{conditions}.

\item [iii)] The 2-cells $\theta:(A, \psi^{\A})\longrightarrow (B, \psi^{\B})$ are just natural transformations $\theta:A\longrightarrow B:\mathcal{C}\longrightarrow\mathcal{D}$ such that the following equation takes place

\begin{equation*}
\psi^{\B}\circ T\theta = \theta S\circ \psi^{\A}
\end{equation*}

\end{itemize}


\subsection{The 2-functor $\Phi^{\mathfrak{\sfM}}_{\E}$}


The 2-functor $\Phi^{\mathfrak{\sfM}}_{\E}:\Adj_{\R}(\mathfrak{C})\longrightarrow \Mnd(\mathfrak{C})$ is defined on $n$-cells as follows

\begin{itemize}

\item [i)] For the 0-cell, $L\dashv R:\mathcal{C}\longrightarrow \mathcal{X}$, $\Phi^{\sfM}_{\E}(L\dashv R) = (\mathcal{C}, RL, R\varepsilon^{\sL\R}L, \eta^{\R\sL})$. That is to say, the monad induced on the domain of the adjunction.

\item [ii)] For the 1-cell, $(J, V, \lambda^{\J\V}):L\dashv R\longrightarrow \overline{L}\dashv \overline{R}$,

\begin{equation*}
\Phi^{\sfM}_{\E}(J, V, \lambda^{\J\V}) = (J, \varrho^{\J\V}\!L\circ\overline{R}\lambda^{\J\V}):(\mathcal{C}, RL)\longrightarrow(\mathcal{D}, \overline{RL})
\end{equation*}

\noindent it will be useful to provide the following notation, $\Phi(\lambda^{\J\V}) = \varrho^{\J\V}\!L\circ\overline{R}\lambda^{\J\V}$.

\item [iii)] For the 2-cell, $(\alpha,\beta)$, $\Phi^{\sfM}_{\E}(\alpha,\beta) = \alpha$.

\end{itemize}


\subsection{The 2-functor $\Psi^{\sfM}_{\E}$}


The 2-functor $\Psi^{\sfM}_{\E}: \Mnd(\mathfrak{C})\longrightarrow \Adj_{\R}(\mathfrak{C})$ can be constructed if the initial 2-category \emph{admits the construction of algebras} \cite{stro_fotm}, which is certainly the case for $\Cat$. It is defined on $n$-cells as follows.

\begin{itemize}

\item [i)] For $(\mathcal{C}, S)$, $\Psi^{\sfM}_{\E}(\mathcal{C}, S) = F^{\sS}\dashv U^{\sS}:\mathcal{C}\longrightarrow\mathcal{C}^{\sS}$. The category $\mathcal{C}^{\sS}$ is the Eilenberg-Moore category for the monad $S$, on $\mathcal{C}$, and the adjunction is the usual free-forgetful adjunction.

\item [ii)] For $(B,\psi^{\B}):(\mathcal{C}, S)\longrightarrow (\mathcal{D}, T)$. \\
$\Psi^{\sfM}_{\E}(B, \psi^{\B}) = (B,\widehat{B}, \lambda^{\B\widehat{\B}}\hspace{0.5pt}):F^{\sS}\dashv U^{\sS}\longrightarrow F^{\T}\dashv U^{\T}$ as in

\begin{equation*}
\begin{array}{ccc}
\xy<1cm,1cm>
\POS (0,0) *+{\mathcal{C}} = "a11",
\POS (30,0) *+{\mathcal{D}} = "a12",
\POS (0,-30) *+{\mathcal{C}^{\sS}} = "a21",
\POS (30,-30) *+{\mathcal{D}^{\T}} = "a22",
\POS (28,-19) *+{} = "x1",
\POS (15,-29) *+{} = "x2",
\POS "a11" \ar^{B} "a12",
\POS "a12" \ar@<-3pt>_{F^{\T}} "a22",
\POS "a22" \ar@<-3pt>_{U^{\T}} "a12",
\POS "a11" \ar@<-3pt>_{F^{\sS}} "a21",
\POS "a21" \ar@<-3pt>_{U^{\sS}} "a11",
\POS "a21" \ar_{\widehat{B}} "a22",
\endxy 
\end{array}
\end{equation*}

\noindent where the functor $\widehat{B}:\mathcal{C}^{\sS}\longrightarrow\mathcal{D}^{\T}$ is defined as

\begin{equation}\label{1703021608}
\widehat{B}(M, k_{\M}) = \big(BU^{\sS}(M, k_{\M}), BU^{\sS}\varepsilon^{\U\!\F\!\sS}(M,k_{\M})\cdot \psi^{\B}U^{\sS}(M,k_{\M})\big)
\end{equation}

$\phantom{br}$

\noindent for $(M,k_{\M})$ in  $\mathcal{C}^{\sS}$. On morphisms, $\widehat{B}(\overline{m}) = \overline{Bm}$. The left mate $\lambda^{\B\widehat{B}}$ fulfills the following equation $U^{\T}\lambda^{\B\widehat{B}} = \psi^{\B}$. It will be useful to make the following notation $\Psi(\psi^{\B}) = \lambda^{\B\widehat{B}}$\\

Another notation for such a functor would be $\mathcal{L}_{\psi}(B) := \widehat{B}$, where the author is considering that $\mathcal{L}$ stands for \emph{lifting}. Also, the notation $B^{\psi}$ is particularly useful. The author will use any of these notations that suits better to the context at hand.\\

The bar over the morphism $m$ means that, although $m$ is in $\mathcal{C}$, it fulfills an additional requirement for algebras. For example, in $U^{\sS}(\overline{m}) = m$, this requirement is forgotten.

\item [iii)] For $\theta:(A,\psi^{\A})\longrightarrow(B,\psi^{\B})$, $\Phi^{\sfM}_{\E}(\theta) = (\theta, \widehat{\theta}\hspace{0.5pt})$, where $\widehat{\theta}: \widehat{A}\longrightarrow\widehat{B}:\mathcal{C}^{\sS}\longrightarrow\mathcal{D}^{\T}$ and whose component, at $(M,k_{\M})$ in $\mathcal{C}^{\sS}$, is

\begin{equation*}
\widehat{\theta}(M, k_{\M}) = \overline{\theta M}
\end{equation*}

\end{itemize}


\subsection{The 2-adjunction $\Phi^{\sfM}_{\E}\dashv\Psi^{\sfM}_{\E}$}


The 2-adjunction can be completed, along the previous pair of 2-functors, by giving the following unit and counit.

\begin{itemize}

\item [i)] The component of the unit $\eta^{\Psi\Phi}: 1_{\Adj_{\R}(\mathfrak{C})}\longrightarrow \Psi^{\sfM}_{\E}\Phi^{\sfM}_{\E}$, on the 0-cell $L\dashv R$, is given by

\begin{equation*}
\eta^{\Psi\Phi}(L\dashv R):= (1_{\C}, \mathfrak{K}_{\mathcal{\D}}^{\R\sL}) = L\dashv R\longrightarrow F^{\R\sL}\dashv U^{\R\sL}
\end{equation*}

\noindent where $\mathfrak{K}_{\mathcal{\D}}^{\R\sL}:\mathcal{D}\longrightarrow \mathcal{C}^{\R\sL}$ is the comparison functor. The notation for this functor tries to codify as much its domain as its codomain, in order to minimize possible confusions.

\item [ii)] The component of the counit $\varepsilon^{\Phi\Psi}: \Phi^{\sfM}_{\E}\Psi^{\sfM}_{\E}\longrightarrow 1_{\Mnd(\mathfrak{C})}$, on the 0-cell $(\mathcal{C}, S)$ is

\begin{equation*}
\varepsilon^{\Phi\Psi}(\mathcal{C}, S):= (1_{\C}, 1_{\sS}): (\mathcal{C}, S)\longrightarrow (\mathcal{C}, S)
\end{equation*}  

\end{itemize}

$\phantom{br}$

All of the previous data make the following Proposition.

\newtheorem{1708231751}{Proposition}[subsection]
\begin{1708231751}
  
There exists a 2-adjunction

\begin{equation*}
 \xy<1cm,1cm>
  \POS (0,0) *+{\Adj_{\R}(\mathfrak{C})} = "a11",
  \POS (35,0) *+{\Mnd(\mathfrak{C})} = "a12",
  \POS "a11" \ar@<-3pt>_{\Phi_{\E}^{\sfM}} "a12",
  \POS "a12" \ar@<-3pt>_{\Psi^{\sfM}_{\E}} "a11"
 \endxy
\end{equation*}

\begin{flushright}
$\Box$
\end{flushright}
  
\end{1708231751}


\section{Adjoint Objects in 2-categories}


In this section, the definitions of \emph{adjoint objects} are developed as much in $\Adj_{\R}(\mathfrak{C})$ as in $\Mnd\:(\mathfrak{C})$.

\subsection{Adjoint Objects in $\Adj_{\R}(\mathfrak{C})$}

In this subsection, a characterization of adjoint objects in the 2-category $\Adj_{\R}(\mathfrak{C})$ is given. These adjoint objects corresponds to the usual definition of an \emph{adjoint object} in a general 2-category $\mathcal{A}$, nevertheless the definition is reviewed in order to characterize these structures.

\newtheorem{1612201652}{Definition}[subsection]
\begin{1612201652}
  
An \emph{adjoint object} in $\Adj_{\R}(\mathfrak{C})$ is comprised of the following.

\begin{enumerate}

\item [i)] A pair of 1-cells

\begin{eqnarray*}
(J,V,\lambda^{\J\V}) &:& L\dashv R\longrightarrow \overline{L}\dashv\overline{R}\:,\\
(K,W,\lambda^{\K\W}) &:& \overline{L}\dashv\overline{R}\longrightarrow L\dashv R\:.
\end{eqnarray*}

\item [ii)] A pair of 2-cells called \emph{unit} and \emph{counit}, respectively

\begin{eqnarray*}
(\eta^{\K\!\J},\eta^{\W\V}) &:& (1_{\mathcal{C}},1_{\mathcal{D}}) \longrightarrow (KJ, WV)\\
(\varepsilon^{\J\!\K},\varepsilon^{\V\W}) &:& (JK, VW)\longrightarrow (1_{\mathcal{X}}, 1_{\mathcal{Y}})
\end{eqnarray*}

\noindent such that they fulfill the following triangular identities

\begin{eqnarray*}
(\varepsilon^{\J\!\K}J\circ J\eta^{\K\J}, \:\varepsilon^{\V\W}V\circ V\eta^{\W\V}) &=& (1_{J}, 1_{V})\\
(K\varepsilon^{\J\!\K}\circ \eta^{\K\J}K, \: W\varepsilon^{\V\W}\circ \eta^{\W\V}W) &=& (1_{K}, 1_{W})  
\end{eqnarray*}

\end{enumerate}

\end{1612201652}

Similar to Theorem 3.13 in \cite{bral_homm}, this type of adjoint object can be characterized by the existence of a natural transformation inverse.\\

\newtheorem{1612211654}[1612201652]{Proposition}
\begin{1612211654}\label{1612211654}

Consider the following diagram in $\Cat$, where $L\dashv R$ and $\overline{L}\dashv\overline{R}$,

\begin{equation*}
\xy<1cm,1cm>
\POS (0,0) *+{\mathcal{C}} = "a11",
\POS (30,0) *+{\mathcal{D}} = "a12",
\POS (0,-30) *+{\mathcal{X}} = "a21",
\POS (30,-30) *+{\mathcal{Y}} = "a22",
\POS "a11" \ar@<-3pt>_{J} "a12",
\POS "a12" \ar@<-3pt>_{K} "a11",
\POS "a12" \ar@<-3pt>_{\overline{L}} "a22",
\POS "a22" \ar@<-3pt>_{\overline{R}} "a12",
\POS "a11" \ar@<-3pt>_{L} "a21",
\POS "a21" \ar@<-3pt>_{R} "a11",
\POS "a21" \ar@<-3pt>_{V} "a22",
\POS "a22" \ar@<-3pt>_{W} "a21",
\endxy
\end{equation*}

$\phantom{br}$

\noindent Consider $J\dashv K$, $V\dashv W$ as classical adjunctions  and $(J,V,\lambda^{\J\V})$ a morphism in $\Adj_{\R}(\mathfrak{C})$. The following assertions are equivalent:

\begin{enumerate}
    
\item [i)] Exists an adjoint object in $\Adj_{\R}(\mathfrak{C})$, where $(K, W)$ is extended to a 1-cell $(K,W,\lambda^{\K\W})$, 

\begin{equation*}
(J,V,\lambda^{\J\V})\dashv(K,W,\lambda^{\K\W})\hspace{1pt}.
\end{equation*}

\item [ii)] $\lambda^{\J\V}$ is invertible.

\end{enumerate}

\noindent In such a case, $\lambda^{\K\W}$ is the mate of the inverse of $\lambda^{\J\V}$. The natural transformation $\lambda^{\K\W}$ might be called \emph{adjoint} of $\lambda^{\J\V}$, the corresponding notation is $\lambda^{\K\W} = ad(\lambda^{\J\V}\hspace{-1pt})$

\end{1612211654}

\noindent \emph{Proof:}\\

\noindent $i \Rightarrow ii$.\\

The proposed inverse, for $\lambda^{\J\V}$, is the following

\begin{equation}\label{1611071526}
\gamma^{\J\V} = \varepsilon^{\V\W}\overline{L}J\circ V\lambda^{\K\W}J\circ VL\eta^{\K\J}
\end{equation}

For example,

\begin{eqnarray*}
\gamma^{\J\V}\circ\lambda^{\J\V} &=& \varepsilon^{\V\W}\overline{L}J\circ V\lambda^{\K\W}J\circ VL\eta^{\K\J}\circ \lambda^{\J\V} \\
&=& \varepsilon^{\V\W}\overline{L}J\circ V\lambda^{\K\W}J\circ \lambda^{\J\V}KJ \circ \overline{L}J\eta^{\K\J}\\
&=& \overline{L}\varepsilon^{\J\K}J\circ \overline{L} J\eta^{\K\J} = \overline{L}(\varepsilon^{\J\K}J\circ J\eta^{\K\J})\\
&=& \overline{L}1_{\J} = 1_{\overline{\sL}\J}
\end{eqnarray*}

$\phantom{br}$

In the third equality, it was used the fact that $(\varepsilon^{\J\K}\!,\:\varepsilon^{\V\W})$ is a 2-cell in $\Adj_{\R}(\mathfrak{C})$. In the fifth one, the triangular identity associated to $J$, for the adjunction $J\dashv K$, was applied.\\

In a similar way, $\lambda^{\J\V}\circ\gamma^{\J\V} = 1_{\V\sL}$.\\

\noindent $ii \Rightarrow i$.\\

Supposing the existence of the inverse, the natural transformation $\lambda^{\K\W}$ is defined as follows

\begin{equation*}
\lambda^{\K\W} = W\overline{L}\varepsilon^{\J\K}\circ W\gamma^{\J\V}\! K\circ \eta^{\W\V}\!LK
\end{equation*}

$\phantom{br}$

That is to say, $\lambda^{\K\W} ={}_{\K}\hspace{-0.3pt}\mathbf{m}_{\W}\hspace{0.6pt}(\:\gamma^{\J\V})$.\\

In order for $(K, W, \lambda^{\K\W})$ to be a morphism in $\Adj_{\R}(\mathfrak{C})$, the mate of $\lambda^{\K\W}$ must be an isomorphic natural transformation. Then, consider the mate of $\lambda^{\K\W}$

\begin{equation*}
\rho^{\K\W} := RW\varepsilon^{\overline{\sL\R}}\circ RW\overline{L}\varepsilon^{\J\K}\overline{R}\circ RW\gamma^{\J\V}\! K\overline{R}\circ R\eta^{\W\V}LK\overline{R}\circ \eta^{\R\sL}K\overline{R}
\end{equation*}

\noindent and its proposed inverse is

\begin{equation*}
\delta^{\K\W} := K\overline{R}\varepsilon^{\V\W} \circ K\overline{R}V \varepsilon^{\sL\R} W\circ K\overline{R}\lambda^{\J\V}RW\circ K\eta^{\overline{\R\sL}}JRW\circ \eta^{\K\!\J}RW
\end{equation*}

$\phantom{br}$

The equation $\rho^{\K\W}\circ\delta^{\K\W} = 1_{\R\W}$ is proved as follows.

\begin{eqnarray*}
&& \rho^{\K\W}\circ\delta^{\K\W} =\\
  &&RW\varepsilon^{\overline{\sL\R}}\circ RW\overline{L}\varepsilon^{\J\K}\overline{R}\circ RW\gamma^{\J\V}K\overline{R}\circ R\eta^{\W\V}LK\overline{R}\circ \eta^{\R\sL}K\overline{R}\circ K\overline{R}\varepsilon^{\V\W} \circ K\overline{R}V \varepsilon^{\sL\R} W\circ\\
  && \qquad\quad\:\circ K\overline{R}\lambda^{\J\V}RW\circ K\eta^{\overline{\R\sL}}JRW\circ \eta^{\K\!\J}RW \\
  && = RW\varepsilon^{\V\W}\circ RWV\varepsilon^{\sL\R}W\circ RW\lambda^{\J\V}RW\circ RW\overline{L}\varepsilon^{\J\K} \overline{R}\overline{L}JRW\circ RW\varepsilon^{\overline{\sL\R}}\:\overline{L}JRW\\
  &&  \qquad\quad RW\overline{L}JK\eta^{\overline{\R\sL}}JRW\circ RW\overline{L}J\eta^{\K\!\J} RW\circ RW \gamma^{\J\V}RW\circ R\eta^{\W\V}L RW\circ\eta^{\R\sL}RW\\
  &&  = \eta^{\R\sL}RW\circ R\eta^{\W\V}LRW\circ RW\gamma^{\J\V}RW\circ RW\lambda^{\J\V} RW\circ RWV\varepsilon^{\sL\R}W\circ RW\varepsilon^{\V\W} \\
  && = \eta^{\R\sL}RW\circ R\eta^{\W\V}LRW\circ RWV\varepsilon^{\sL\R}W\circ RW\varepsilon^{\V\W}\\
  && = 1_{\R\W}
\end{eqnarray*}

In the third equality, the triangular identity associated to $\overline{L}J$ of the composed adjunction $\overline{L}J\dashv K\overline{R}$ was used. In the fifth, the triangular identity associated to $RW$ was applied.\\

In a similar fashion, it can be proved that $\delta^{\K\W}\circ\rho^{\K\W} = 1_{\K\overline{\R}}$. Therefore, $(K,W,\lambda^{\K\W})$ is a morphism, or 1-cell, in $\Adj_{\R}(\mathfrak{C})$.\\

Remains to prove that the pair $(\eta^{\K\!\J}\!,\:\eta^{\W\V}):(1_{\mathcal{C}}, 1_{\mathcal{X}})\longrightarrow (KJ, WV, W\lambda^{\J\V}\circ\lambda^{\K\W}J):L\dashv R\longrightarrow L\dashv R$ is a 2-cell in $\Adj_{\R}(\mathfrak{C})$. In particular, it is required that

\begin{equation*}
W\lambda^{\J\V}\circ\lambda^{\K\W}J\circ L\eta^{\K\!\J} = \eta^{\W\V}L
\end{equation*}

Therefore,

\begin{eqnarray*}
 W\lambda^{\J\V}\circ\lambda^{\K\W}J\circ L\eta^{\K\J} &=& W\lambda^{\J\V}\circ W\overline{L}\varepsilon^{\J\K}J\circ W\gamma^{\J\V}KJ\circ \eta^{\W\V}LKJ\circ L\eta^{\K\!\J} \\
&=& W\lambda^{\J\V}\circ W\overline{L}(\varepsilon^{\J\!\K}J\circ J\eta^{\K\!\J})\circ W\gamma^{\J\V}\circ\eta^{\W\V}L\\
&=& W\lambda^{\J\V}\circ W\gamma^{\J\V}\circ \eta^{\W\V}L = \eta^{\W\V}L
\end{eqnarray*}

In the third equality, it was used the triangular identity, of $J\dashv K$ , associated to $J$.\\

That the pair $(\varepsilon^{\J\K}\!, \:\varepsilon^{\V\W}):(JK, VW,V\lambda^{\K\W}\circ\lambda^{\J\V}K):L\dashv R\longrightarrow L\dashv R$ is a 2-cell in $\Adj_{\R}(\mathfrak{C})$ is proved similarly. Finally, the triangular identities are fulfilled since composition, and whiskering, of 2-cells in $\Adj_{R}(\mathfrak{C})$ are composed, and whiskered, componently as in $\Cat$.

\begin{flushright}
$\Box$
\end{flushright}


\subsection{Adjoint Objects in $\Mnd(\mathfrak{C})$}


As in the previous section, a detailed account of \emph{adjoint objects} in the 2-category $\Mnd(\mathfrak{C})$ is given.\\

\newtheorem{1612212232}{Definition}[subsection]
\begin{1612212232}

An \emph{adjoint object} in $\Mnd(\mathfrak{C})$ is comprised of the following items:

\begin{enumerate}

\item [i)] A pair of 1-cells,

\begin{eqnarray*}
(J, \psi^{\J}) &:& (\mathcal{C}, S) \longrightarrow (\mathcal{D}, T)\:,\\
(K, \psi^{\K}) &:& (\mathcal{D}, T) \longrightarrow (\mathcal{C}, S)\:.
\end{eqnarray*}

\item [ii)] A pair of 2-cells, the \emph{unit} and the \emph{counit} of the adjoint object

\begin{eqnarray*}
\eta^{\K\!\J}&:&(1_{\mathcal{C}}, 1_{S})\longrightarrow (KJ, K\psi^{\J}\circ\psi^{\K}\!J):(\mathcal{C}, S)\longrightarrow(\mathcal{C}, S)\\
\varepsilon^{\J\!\K} &:& (JK, J\psi^{\K}\circ\psi^{\J}K) \longrightarrow (1_{\mathcal{D}}, 1_{\T}):(\mathcal{D}, T)\longrightarrow (\mathcal{D}, T) 
\end{eqnarray*}

\noindent such that they fulfill the triangular identities

\begin{eqnarray*}
\varepsilon^{\J\!\K}J\circ J\eta^{\K\!\J}&=& 1_{\J}\\
K\varepsilon^{\J\!\K}\circ\eta^{\K\!\J}K &=& 1_{\K}
\end{eqnarray*}

\end{enumerate}

\end{1612212232}

This type of object can be characterised using the Theorem 3.13 in \cite{bral_homm}. However, it is restated and proved again within the context of this article.\\

\newtheorem{1612221117}[1612212232]{Proposition}
\begin{1612221117}\label{1612221117}
  
Consider the following adjunction $J\dashv K:\mathcal{C}\longrightarrow\mathcal{D}$, such that $J$ is part of the 1-cell

\begin{equation}\label{1611091846}
\xy<1cm,1cm>
\POS (0,0) *+{(\mathcal{C}, S)} = "a11",
\POS (30,0) *+{(\mathcal{D}, T)} = "a12",
\POS "a11" \ar^{(J,\: \psi^{\J})} "a12",
\endxy
\end{equation}

$\phantom{br}$

\noindent in $\Mnd(\mathfrak{C})$. Then the following assertions are equivalent:

\begin{enumerate}

\item [i)] Exists an adjoint object in $\Mnd(\mathfrak{C})$, where $K$ is extended to a 1-cell $(K, \psi^{\K})$,

\begin{equation*}
(J,\psi^{\J})\dashv(K,\psi^{\K})
\end{equation*}

\item [ii)] $\psi^{\J}$ is invertible.

\end{enumerate}

\noindent in such case, $\psi^{\K} = ad(\psi^{\J})$.
  
\end{1612221117}

$\phantom{br}$

\noindent \emph{Proof:}\\

\noindent $i \Rightarrow ii$.\\

The definition of the inverse $\zeta^{\J}$, of $\psi^{\J}$, goes as follows

\begin{equation*}
\zeta^{\J} := \varepsilon^{\J\!\K}TJ\circ J\psi^{\K}J\circ JS\eta^{\K\!\J}
\end{equation*}\\

The equality $\psi^{\J}\circ\zeta^{\J} = 1_{\J\sS}$ is proved

\begin{eqnarray*}
  \psi^{\J}\circ\zeta^{\J} &=& \psi^{\J}\circ \varepsilon^{\J\!\K}TJ\circ J\psi^{\K}J\circ JS\eta^{\K\!\J}\\
  &=& \varepsilon^{\J\!\K}JS\circ JK\psi^{\J}\circ J\psi^{\K}J\circ JS\eta^{\K\!\J}\\
  &=& \varepsilon^{\J\!\K}JS\circ J\eta^{\K\!\J}S \\
  &=& 1_{\J\sS}
\end{eqnarray*}

$\phantom{br}$

\noindent In the third equality, it was used the fact that $\eta^{\K\!\J}$ is a 2-cell in the 2-category $\Mnd(\mathfrak{C})$. In the fourth one, the triangular identity associated to $J$. The case $\zeta^{\J}\circ\psi^{\J} = 1_{\T\!\J}$ is done similarly.\\

\noindent $ii \Rightarrow i$.\\

Suppose that $\psi^{\J}$ has an inverse, $\zeta^{\J}$, then $\psi^{\K}$ is defined as follows

\begin{equation*}
\psi^{\K} := KT\varepsilon^{\J\!\K}\circ K\zeta^{\J}K\circ \eta^{\K\!\J}SK
\end{equation*}\\

First, it is proved that $(K,\psi^{\K}):(\mathcal{D},T)\longrightarrow (\mathcal{C},S)$ is a morphism in $\Mnd(\mathfrak{C})$. In order to do so, the compatibility with the products, \emph{i.e.} $\psi^{\K}\circ\mu^{\sS}K = K\mu^{\T}\circ \psi^{\K}T\circ S\psi^{\K}$, is checked.

\begin{eqnarray*}
 \psi^{\K}\circ\mu^{\sS}K &=& KT\varepsilon^{\J\!\K}\circ K\zeta^{\J}K\circ\eta^{\K\!\J}SK\circ \mu^{\sS}K\\
  &=& KT\varepsilon^{\J\!\K}\circ K\mu^{\T}JK\circ KT\zeta^{\J}K\circ K\zeta^{\J}SK\circ \eta^{\K\!\J}SSK\\
  &=& K\mu^{\T}\circ KTT\varepsilon^{\J\!\K}\circ KT\zeta^{\J}K\circ KT\varepsilon^{\K\!\J}JSK\circ KTJ\eta^{\K\!\J}SK\circ K\zeta^{\J}SK\circ \eta^{\K\!\J}SSK\\
  &=& K\mu^{\T}\circ KT\varepsilon^{\K\!\J}T\circ K\zeta^{\J}KT\circ \eta^{\K\!\J}SKT\circ SKT\varepsilon^{\J\!\K}\circ SK\zeta^{\J}K\circ S\eta^{\K\!\J}SK\\
  &=& K\mu^{\T}\circ\psi^{\K}T\circ S\psi^{\K}.\\
\end{eqnarray*}

If $(J, \psi^{\J}):(\mathcal{C},S)\longrightarrow(\mathcal{D},T)$ is a morphism in $\Mnd(\mathfrak{C})$, then $(J,\zeta^{\J}): (\mathcal{C},S)\longrightarrow(\mathcal{D},T)$ is a morphism in the Kleisli dual of $\Mnd(\mathfrak{C})$. In particular, $\zeta^{\J}\circ J\mu^{\sS} = \mu^{\T}J\circ T\zeta^{\J}\circ \zeta^{\J}S$. This compatibility was used for the second equality. In the fifth equality, it was used the $J$ triangular identity and the definition of $\psi^{\K}$. \\

The compatibility of the units is left to the reader.\\

Next, it is needed that $\eta^{\K\!\J}\! :(1_{\C}, 1_{\sS})\longrightarrow (KJ, K\psi^{\J}\circ\psi^{\K}J):(\mathcal{C},S)\longrightarrow(\mathcal{C},S)$, be a 2-cell in $\Mnd(\mathfrak{C})$, \emph{i.e.} the following equality takes place $K\psi^{\J}\circ\psi^{\K}J\circ S\eta^{\K\!\J} = \eta^{\K\!\J}S$.

\begin{eqnarray*}
 K\psi^{\J}\circ\psi^{\K}J\circ S\eta^{\K\!\J} &=& K\psi^{\J}\circ KT\varepsilon^{\J\!\K}J\circ K\zeta^{\J}KJ\circ \eta^{\K\!\J}SKJ\circ S\eta^{\K\!\J} \\
  &=& K\psi^{\J}\circ KT\varepsilon^{\J\!\K}J\circ KTJ\eta^{\K\!\J}\circ K\zeta^{\J}\circ\eta^{\K\!\J}S\\
  &=& K\psi^{\J}\circ K\zeta^{\J}\circ\eta^{\K\!\J}S\\
  &=& \eta^{\K\!\J}S
\end{eqnarray*}

\noindent In the third equality, it was used the triangular identity associated to $J$. In the fourth one, the fact that $\zeta^{\J}$ is the inverse of $\psi^{\J}$ was applied.\\

Likewise, $\varepsilon^{\J\!\K}:(JK, J\psi^{\K}\circ\psi^{\J}K)\longrightarrow(1_{\smD},1_{\T}):(\mathcal{D},T)\longrightarrow(\mathcal{D},T)$ is a 2-cell in $\Mnd\hspace{0.5pt}(\mathfrak{C})$. Since the composition of 2-cells in $\Mnd\hspace{0.5pt}(\mathfrak{C})$, and the whiskering, is done as in the subjacent 2-category $\Cat$, then the triangular identities are fulfilled.
  
\begin{flushright}
$\Box$
\end{flushright}

In Example 3.12 in \cite{bral_homm}, the \emph{lift of an adjunction} corresponds to an adjoint object in $\Mnd(\mathfrak{C})$. For example, conditions 3a-3d correspond to $G$ and $V$, along with $\zeta$ and $\xi$, being morphisms in $\Mnd(\mathfrak{C})$ and 3e-3f are the requirements for the unit and counit $(h,e)$ being 2-cells in the same 2-category.\\

The results on adjoint objects, using the 2-adjunction $\Phi_{\E}^{\sfM} \dashv \Psi_{\E}^{\sfM}$, can be combined. Take the 0-cells $F^{\sS}\dashv U^{\sS}$ in $\Adj_{\R}(\mathfrak{C})$ and $(\mathcal{D}, T)$ in $\Mnd\hspace{0.5pt}(\mathfrak{C})$. Therefore, there exists an isomorphism of categories, natural in $F^{\sS}\dashv U^{\sS}$ and $(\mathcal{D}, T)$

\begin{equation*}
Hom_{\Adj_{\R}(\mathfrak{C})}\big(F^{\sS}\dashv U^{\sS}, \Psi^{\mathfrak{M}}_{\E}(\mathcal{D}, T)\big) \cong Hom_{\Mnd(\mathfrak{C})}(\Phi^{\mathfrak{M}}_{\E}(F^{\sS}\dashv U^{\sS}), (\mathcal{D}, T))
\end{equation*}

$\phantom{br}$

\noindent in particular, there is a bijection between adjoint objects inside each category. If we take into account the proofs of Proposition \ref{1612211654} and Proposition \ref{1612221117} and the previous isomorphism of categories, we are left, without any need of a proof, with the following Theorem.


\newtheorem{1612221241}[1612212232]{Theorem}
\begin{1612221241}

The following statements are equivalent

\begin{itemize}

\item [i)] There exist adjoint objects in $\Mnd(\mathfrak{C})$ of the form

\begin{equation*}
\xy<1cm,1cm>
\POS (0,0) *+{(\mathcal{C}, S)} = "a11",
\POS (35,0) *+{(\mathcal{D}, T)} = "a12",
\POS "a11" \ar@<-3pt>_{(J, \:\psi^{\J})} "a12",
\POS "a12" \ar@<-3pt>_{(K, \:\psi^{\K})} "a11",
\endxy
\end{equation*}

\item [ii)] There exist adjoint objects in $\Adj_{\R}(\mathfrak{C})$ of the form

\begin{equation*}
\xy<1cm,1cm>
\POS (0,0) *+{\mathcal{C}} = "a11",
\POS (30,0) *+{\mathcal{D}} = "a12",
\POS (0,-30) *+{\mathcal{C}^{\sS}} = "a21",
\POS (30,-30) *+{\mathcal{D}^{\T}} = "a22",
\POS "a11" \ar@<-3pt>_{J} "a12",
\POS "a12" \ar@<-3pt>_{K} "a11",
\POS "a12" \ar@<-3pt>_{F^{\T}} "a22",
\POS "a22" \ar@<-3pt>_{U^{\T}} "a12",
\POS "a11" \ar@<-3pt>_{F^{\sS}} "a21",
\POS "a21" \ar@<-3pt>_{U^{\sS}} "a11",
\POS "a21" \ar@<-3pt>_{\widehat{J}} "a22",
\POS "a22" \ar@<-3pt>_{\widehat{K}} "a21",
\endxy
\end{equation*}

\item [iii)] The natural transformation $\psi^{\J}: TJ\longrightarrow JS$ is invertible.

\item [iv)] The natural transformation $\lambda^{\J\:\widehat{\!\!\J}}:F^{\T}\!J\longrightarrow \widehat{J}F^{\sS}$ is invertible.

\end{itemize}

\end{1612221241}

\begin{flushright}
$\Box$
\end{flushright}


\section{Left Hopf 1-cells and Parametric Adjoint Objects in $\Adj_{\R}(\mathfrak{C})$}


Recalling that the main objective in this article is to characterize parametric adjoint objects as much in $\Adj_{\R}(\mathfrak{C})$ as in $\Mnd(\mathfrak{C})$ and relate them through the 2-adjunction $\Phi_{\E}^{\sfM}\dashv\Psi_{\E}^{\sfM}$, therefore the definition and characterization of these structures in $\Cat$ has to be given. 

\subsection{Preliminars}

The definition of a parametric adjunction is recalled along with the corresponding theorem that characterizes it, \cite{masa_cawm_2ned}.

\newtheorem{1612221311}{Definition}[subsection]
\begin{1612221311}\label{1612221311}

Consider the following categories $\mathcal{P}$, $\mathcal{C}$ and $\mathcal{D}$. A \emph{parametric adjunction}, by $\mathcal{P}$, is a pair of functors of the form

\begin{eqnarray*}
F&:&\mathcal{C}\times\mathcal{P}\longrightarrow\mathcal{D},\\
G&:&\mathcal{P}^{op}\times\mathcal{D}\longrightarrow\mathcal{C},
\end{eqnarray*}

\noindent such that for any $P$ in $\mathcal{P}$, there is an adjunction $F_{\sP}\dashv G_{\sP}$ and, for $p:P\longrightarrow Q$, a \emph{conjugate morphism} of adjunctions $\widetilde{p}:F_{\sP}\dashv G_{\sP}\longrightarrow F_{\Q}\dashv G_{\Q}$. This parametric adjunction can be denoted as $F\dashv_{\:\mathcal{P}} G:\mathcal{C}\longrightarrow\mathcal{D}$.
  
\end{1612221311}

Now, the corresponding characterizing theorem.

\newtheorem{1709020912}[1612221311]{Theorem}
\begin{1709020912}\label{1709030912}
  
Consider a functor $F:\mathcal{C}\times\mathcal{P}\longrightarrow\mathcal{D}$ such that for every $P$ in $\mathcal{P}$ there exists  a functor $G_{\sP}:\mathcal{D}\longrightarrow\mathcal{C}$ and an adjunction 

\begin{equation*}
\xy<1cm,1cm>
\POS (0,0) *+{\mathcal{C}} = "a11",
\POS (28,0) *+{\mathcal{D}} = "a21",
\POS "a11" \ar@<-3pt>_{F_{\sP}} "a21",
\POS "a21" \ar@<-3pt>_{G_{\sP}} "a11",
\endxy
\end{equation*}

\noindent Therefore, exists a unique $G:\mathcal{P}^{op}\times\mathcal{D}\longrightarrow\mathcal{C}$ such that for $P$ 

\begin{equation*}
  G(P,\sim) := G_{\sP}
\end{equation*}

\noindent and for $p^{op}:P'\longrightarrow P$, in $\mathcal{P}^{op}$, a natural transformation

\begin{equation*}  
G(p^{op},\sim) :   G_{\sP'}\longrightarrow G_{\sP},
\end{equation*}

$\phantom{br}$

\noindent further denoted as $G_{p^{op}}$, such that

\begin{equation*}
G_{p^{op}} = G_{\sP}\hspace{0.5pt}\varepsilon^{\F\G\sP'} \circ G_{\sP}F_{p}G_{\sP'}\circ\eta^{\G\F\!\sP}G_{\sP'}.
\end{equation*}

\begin{flushright}
$\Box$
\end{flushright}

\end{1709020912}

The departure from the parametric adjoint objects in $\Cat$ to the 2-category realm is given by the \emph{comonoidal adjunction}, cf. \cite{moie_motc}, and the Hopf adjunction, cf. \cite{bral_homm}

\newtheorem{1708081819}[1612221311]{Definition}
\begin{1708081819}
  
A \emph{comonoidal adjunction} is defined as an adjunction $L\dashv R:\mathcal{C}\longrightarrow\mathcal{D}$ where $\mathcal{C}$ and $\mathcal{D}$ are monoidal categories and $L$ and $R$ are comonoidal functors and the unit and counit $\eta^{\R\sL}:1_{\mathcal{C}}\longrightarrow RL$ and $\varepsilon^{\sL\R}:LR\longrightarrow 1_{\mathcal{D}}$ are natural comonoidal transformations.
\end{1708081819}

In \cite{lojt_apke} there is a characterization of this comonoidal adjunctions. 

\newtheorem{1708081833}[1612221311]{Proposition}
\begin{1708081833}\label{1708081833}

The following statements are equivalent

\begin{itemize}

\item [i)] The adjunction $L\dashv R:\mathcal{C}\longrightarrow\mathcal{D}$ is comonoidal.
  
\item [ii)] The following diagram

\begin{equation*}
\xy<1cm,1cm>
\POS (0,0) *+{\mathcal{C}\times\mathcal{C}} = "a11",
\POS (30,0) *+{\mathcal{C}} = "a12",
\POS (0,-25) *+{\mathcal{D}\times\mathcal{D}} = "a21",
\POS (30,-25) *+{\mathcal{D}} = "a22",
\POS "a11" \ar@<-3pt>_{L\times L} "a21",
\POS "a21" \ar@<-3pt>_{R\times R} "a11",
\POS "a21" \ar_-{\otimes_{\mathcal{D}}} "a22",
\POS "a11" \ar^-{\otimes_{\mathcal{C}}} "a12",
\POS "a12" \ar@<-3pt>_{L} "a22",
\POS "a22" \ar@<-3pt>_{R} "a12",
\endxy
\end{equation*}

\noindent is a 1-cell, $(\otimes_{\mathcal{C}}, \otimes_{\mathcal{D}}, \lambda^{\otimes\mathcal{C}\mathcal{D}})$   , in $\Adj_{\R}(\mathfrak{C})$.

\end{itemize}

\end{1708081833}

Let us remember the definition of a \emph{Hopf operator} in order to start the extension of this concepts to the context of 2-categories.

\newtheorem{1708081900}[1612221311]{Definition}
\begin{1708081900}

Let $L\dashv R:\mathcal{C}\longrightarrow \mathcal{D}$ be a comonoidal adjunction. The \emph{left Hopf operator}, $\mathfrak{H}$ is the following natural transformation

\begin{equation}
\mathfrak{H}(\lambda^{\otimes\mathcal{C}\mathcal{D}}):\otimes_{\mathcal{D}}(L\times \varepsilon^{\sL\!\R})\circ \lambda^{\otimes\mathcal{C}\mathcal{D}}(\mathcal{C}\times R):\mathcal{C}\times\mathcal{D}\longrightarrow\mathcal{D}
\end{equation}

\begin{equation*}
 \xy<1cm,1cm>
   \POS (0,0) *+{\mathcal{C}\times\mathcal{D}} = "a11",
   \POS (25,0) *+{\mathcal{C}\times\mathcal{C}} = "a12",
   \POS (55,0) *+{\mathcal{C}} = "a13",
   \POS (25,-25) *+{\mathcal{D}\times\mathcal{D}} = "a22",
   \POS (55,-25) *+{\mathcal{D}} = "a23",
   \POS (53,-15) *+{} = "x1",
   \POS (40,-23) *+{} = "x2",
   \POS (25,-5) *+{} = "y1",
   \POS (11,-11) *+{} = "y2",
   \POS "a11" \ar_{L\times\mathcal{D}} "a22",
   \POS "a22" \ar_-{\otimes_{\mathcal{D}}} "a23",
   \POS "a11" \ar^{\mathcal{C}\times R} "a12",
   \POS "a12" \ar_{L\times L} "a22",
   \POS "a12" \ar^-{\otimes_{\mathcal{C}}} "a13",
   \POS "a13" \ar^{L} "a23",
   \POS "x1" \ar@/_0.5pc/_{\lambda^{\otimes\mathcal{C}\mathcal{D}}} "x2",
   \POS "y1" \ar@/_0.4pc/_{L\times\varepsilon^{\sL\!\R}} "y2",
 \endxy
\end{equation*}

\end{1708081900}


\subsection{Hopf 1-cells}


The objective of this section is to extend the definition of a parametric adjunction to the 2-category of adjunctions $\Adj_{\R}(\mathfrak{C})$.\\

Consider a 1-cell in $\Adj_{\R}(\mathfrak{C})$ of the form $(J, V, \lambda^{\J\V}):L\times\widetilde{L}\dashv R\times\widetilde{R}\longrightarrow \overline{L}\dashv\overline{R}$

 \begin{equation*}
      \xy<1cm,1cm>
      \POS (0,0) *+{\mathcal{C}\times\mathcal{P}} = "a11",
      \POS (30,0) *+{\mathcal{D}} = "a12",
      \POS (0,-25) *+{\mathcal{X}\times\mathcal{Q}} = "a21",
      \POS (30,-25) *+{\mathcal{Y}} = "a22",
      \POS "a11" \ar@<-3pt>_{L\times \widetilde{L}} "a21",
      \POS "a21" \ar@<-3pt>_{R\times \widetilde{R}} "a11",
      \POS "a21" \ar_-{V} "a22",
      \POS "a11" \ar^-{J} "a12",
      \POS "a12" \ar@<-3pt>_{\overline{L}} "a22",
      \POS "a22" \ar@<-3pt>_{\overline{R}} "a12",
       \endxy
 \end{equation*}\\

Suppose that the functors $J:\mathcal{C}\times\mathcal{P}\longrightarrow\mathcal{D}$ and $V:\mathcal{X}\times\mathcal{Q}\longrightarrow\mathcal{Y}$ are part of classical parametric adjunctions, namely $J\dashv_{\hspace{0.5pt}\mathcal{P}} K$ and $V\dashv_{\hspace{0.5pt}\mathcal{Q}} W$. There is no immediate translation of a parametric adjoint object to the 2-category $ \Adj_{\R}(\mathfrak{C})$ due to a little obstacle. The problem arises with the possible definition of the 1-cell $(K, W, \lambda^{\K\W})$ where the opposite adjunction, for $\widetilde{L}\dashv \widetilde{R}$, $\widetilde{R}^{op}\dashv \widetilde{L}^{op}:\mathcal{Q}^{op}\longrightarrow\mathcal{P}^{op}$ change the domain and the codomain, therefore a 1-cell of the form $(K, W, \lambda^{\K\W})$ cannot be defined. \\

Hence, the objective can be changed to the study of what extension a parametric adjunction can be reasoning within the 2-category $\Adj_{\R}(\mathfrak{C})$. For that, the following modifications of definitions, in \cite{bral_homm}, can be given.
 
\newtheorem{1708081705}{Definition}[subsection]
\begin{1708081705}
  
Let $(J, V, \lambda^{\J\V}):L\times\widetilde{L}\dashv R\times\widetilde{R}\longrightarrow \overline{L}\dashv\overline{R}$ be a 1-cell in $\Adj_{\R}(\mathfrak{C})$. A \emph{left Hopf operator} $\mathfrak{H}$ on $(J, V, \lambda^{\J\V})$ is a 1-cell in $\Adj_{\R}(\mathfrak{C})$ of the form

\begin{equation*}
\mathfrak{H}(J, V, \lambda^{\J\V}):= (J(\mathcal{C}\times\widetilde{R}), V, H(\lambda^{\J\V})): L\times\mathcal{Q}\dashv R\times\mathcal{Q}\longrightarrow \overline{L}\dashv\overline{R}
\end{equation*}

$\phantom{br}$

\noindent where $H(\lambda^{\J\V})$ is the following natural transformation
  
\begin{equation*}
 \xy<1cm,1cm>
  \POS (0,0) *+{\mathcal{C}\times\mathcal{Q}} = "a11",
  \POS (25,0) *+{\mathcal{C}\times\mathcal{P}} = "a12",
  \POS (55,0) *+{\mathcal{D}} = "a13",
  \POS (25,-25) *+{\mathcal{X}\times\mathcal{Q}} = "a22",
  \POS (55,-25) *+{\mathcal{Y}} = "a23",
  \POS (53,-15) *+{} = "x1",
  \POS (40,-23) *+{} = "x2",
  \POS (25,-8) *+{} = "y1",
  \POS (15,-15) *+{} = "y2",
  \POS "a11" \ar_{L\times\mathcal{Q}} "a22",
  \POS "a22" \ar_-{V} "a23",
  \POS "a11" \ar^{\mathcal{C}\times\widetilde{R}} "a12",
  \POS "a12" \ar^{L\times\widetilde{L}} "a22",
  \POS "a12" \ar^-{J} "a13",
  \POS "a13" \ar^{\overline{L}} "a23",
  \POS "x1" \ar@/_0.5pc/_{\lambda^{\J\V}} "x2",
  \POS "y1" \ar@/_0.4pc/_{L\times\varepsilon^{\widetilde{L}\widetilde{R}}} "y2",
\endxy
\end{equation*}

\noindent that is to say

\begin{equation}
H(\lambda^{\J\V}) = (L\times\varepsilon^{\widetilde{\sL}\widetilde{\R}})\cdot_{p}\hspace{0.4pt}\lambda^{\J\V} = V(L\times \varepsilon^{\widetilde{\sL}\widetilde{\R}})\circ \lambda^{\J\V}(\mathcal{C}\times\widetilde{R})\:.
\end{equation}

\end{1708081705}

$\phantom{br}$

\newtheorem{1708081810}[1708081705]{Definition}
\begin{1708081810}

A \emph{left Hopf} 1-cell in $\Adj_{\R}(\mathfrak{C})$,

\begin{equation*}
(J, V, \lambda^{\J\V}):L\times\widetilde{L}\dashv R\times\widetilde{R}\longrightarrow \overline{L}\dashv\overline{R},
\end{equation*}

$\phantom{br}$

\noindent is such that $H(\lambda^{\J\V})$ is invertible. In such a case, the inverse is denoted as $N(\lambda^{\J\V})$. 

\end{1708081810}

Consider a left Hopf 1-cell in $\Adj_{\R}(\mathfrak{C})$, $(J, V, \lambda^{\J\V})$, therefore its left Hopf operator is invertible and so is the following natural transformation, for any $Q$ in $\mathcal{Q}$,

\begin{equation*}
 \xy<1cm,1cm>
  \POS (0,0) *+{\mathcal{C}} = "a11",
  \POS (30,0) *+{\mathcal{C}\times \mathbf{1}} = "a12",
  \POS (60,0) *+{\mathcal{C}\times\mathcal{Q}} = "a13",
  \POS (90,0) *+{\mathcal{C}\times\mathcal{P}} = "a14",
  \POS (120,0) *+{\mathcal{D}} = "a15",
  \POS (0,-25) *+{\mathcal{X}} = "a21",
  \POS (30,-25) *+{\mathcal{X}\times \mathbf{1}} = "a22",
  \POS (60,-25) *+{\mathcal{X}\times\mathcal{Q}} = "a23",
  \POS (90,-25) *+{\mathcal{X}\times\mathcal{Q}} = "a24",
  \POS (88,-15) *+{} = "x1",
  \POS (75,-23) *+{} = "x2",
  \POS (118,-15) *+{} = "y1",
  \POS (105,-23) *+{} = "y2",
  \POS (120,-25) *+{\mathcal{Y}} = "a25",
  \POS "a11" \ar_{L} "a21",
  \POS "a21" \ar_{\rho^{-1}\mathcal{X}} "a22",
  \POS "a22" \ar_{\mathcal{X}\times E_{\Q}} "a23",
  \POS "a23" \ar_{1_{\mathcal{X}}\times 1_{\mathcal{Q}}} "a24",
  \POS "a24" \ar_-{V} "a25",
  \POS "a11" \ar^{\rho^{-1}\mathcal{C}} "a12",
  \POS "a12" \ar^{\mathcal{C}\times E_{\Q}} "a13",
  \POS "a13" \ar^{\mathcal{C}\times \widetilde{R}} "a14",
  \POS "a14" \ar^-{J} "a15",
  \POS "a15" \ar^{\overline{L}} "a25",
  \POS "a12" \ar^{L\times\mathbf{1}} "a22",
  \POS "a13" \ar^{L\times\mathcal{Q}} "a23",
  \POS "a14" \ar^{L\times\widetilde{L}} "a24",
  \POS "a15" \ar^{\overline{L}} "a25",
  \POS "x1" \ar@/_0.5pc/_{L\times\varepsilon^{\widetilde{\sL}\widetilde{\R}}} "x2",
  \POS "y1" \ar@/_0.5pc/_{\lambda^{\J\V}} "y2",
 \endxy
\end{equation*}

\noindent The functor $E_{Q}$ stands for evaluation at $Q$ in $\mathcal{Q}$. The previous natural transformation can be written as follows

\begin{equation*}
 \xy<1cm,1cm>
  \POS (0,0) *+{\mathcal{C}} = "a11",
  \POS (30,0) *+{\mathcal{D}} = "a12",
  \POS (0,-25) *+{\mathcal{X}} = "a21",
  \POS (30,-25) *+{\mathcal{Y}} = "a22",
  \POS (28,-15) *+{} = "x1",
  \POS (15,-23) *+{} = "x2",
  \POS "a11" \ar_{L} "a21",
  \POS "a21" \ar_{V_{\Q}} "a22",
  \POS "a11" \ar^{J_{\widetilde{\R}\Q}} "a12",
  \POS "a12" \ar^{\overline{L}} "a22",
  \POS "x1" \ar@/_0.5pc/_{\lambda^{\J\V\!\widetilde{\R}\Q}} "x2",
\endxy
\end{equation*}

\textbf{Remark}: If $J\dashv_{\hspace{0.5pt}\mathcal{P}} K$ then $J(\mathcal{C}\times\widetilde{R})\dashv_{\hspace{0.4pt}\mathcal{Q}}K(\widetilde{R}^{op}\times\mathcal{D})$.\\

Due to the previous remark, there are two parametric adjunctions on $\mathcal{Q}$, $J(\mathcal{C}\times\widetilde{R})\dashv_{\hspace{0.4pt}\mathcal{Q}}K(\widetilde{R}^{op}\times\mathcal{D})$ and $V \dashv_{\mathcal{Q}} W$ with corresponding adjunctions $J_{\widetilde{R}Q}\dashv K_{\widetilde{R}Q}$ and $V_{Q}\dashv W_{\Q}$, for any $Q$ in $\mathcal{Q}$, such that

\begin{equation*}
(J_{\widetilde{\R}\Q}, V_{\Q}, \lambda^{\J\V\!\widetilde{R}\Q})
\end{equation*}

$\phantom{br}$

\noindent is a 1-cell in $\Adj_{\R}(\mathfrak{C})$ and $\lambda^{\J\V\!\widetilde{R}\Q}$ is invertible. If Proposition \ref{1612211654} is recalled for this situation, there exists an adjoint object in $\Adj_{\R}(\mathfrak{C})$

\begin{equation*}
(J_{\widetilde{\R}\Q}, V_{\Q}, \lambda^{\J\V\!\widetilde{\R}\Q})\dashv(K_{\widetilde{\R}\Q}, W_{\Q}, \lambda^{\K\W\!\widetilde{\R}\Q})
\end{equation*}

$\phantom{br}$

\noindent where $\lambda^{\K\W\!\widetilde{\R}\Q} = ad(\lambda^{\J\V\!\widetilde{\R}\Q})$.\\

The natural transformation $\lambda^{\K\W\!\widetilde{\R}\Q}:LK_{\widetilde{\R}\Q}\longrightarrow W_{\Q}\overline{L}:\mathcal{D}\longrightarrow\mathcal{X}$ can be extended to a dinatural transformation of the form

\begin{equation*}
\lambda^{\K\W\!\widetilde{\R}}: LK(\widetilde{R}^{op}\times\mathcal{D})\longrightarrow W(\mathcal{Q}^{op}\times\overline{L}):\mathcal{Q}^{op}\times\mathcal{D}\longrightarrow\mathcal{X}
\end{equation*}

$\phantom{br}$

This claim is stated as the following proposition\\

\newtheorem{1708091812}[1708081705]{Proposition}
\begin{1708091812}\label{1708091812}

In the previous context, there exists a dinatural transformation of the form

\begin{equation*}
\lambda^{\K\W\!\widetilde{\R}}: LK(\widetilde{R}^{op}\times\mathcal{D})\longrightarrow W(\mathcal{Q}^{op}\times\overline{L}):\mathcal{Q}^{op}\times\mathcal{D}\longrightarrow\mathcal{X}
\end{equation*}

$\phantom{br}$

\noindent defined, on $(Q,D)$ in $\mathcal{Q}^{op}\times\mathcal{D}$, as $\lambda^{\K\W\widetilde{\R}}(Q, D):= \lambda^{\K\W\widetilde{\R}\Q}D$ and such that for any $(q^{op},d):(Q',D)\longrightarrow(Q,D')$ in $\mathcal{Q}^{op}\times\mathcal{D}$, the following diagram commutes

\begin{equation*}
  \xy<1cm, 1cm>
   \POS (25,0) *+{LK(\widetilde{R}Q', D')} = "a12",
   \POS (70,0) *+{W(Q',\overline{L}D')} = "a13",
   \POS (0,-15) *+{LK(\widetilde{R}Q', D)} = "a21",
   \POS (95,-15) *+{W(Q,\overline{L}D')} = "a24",
   \POS (25,-30) *+{LK(\widetilde{R}Q, D)} = "a32",
   \POS (70,-30) *+{W(Q,\overline{L}D)} = "a33",
   \POS "a21" \ar_{LK((\widetilde{R}q)^{op}, D)} "a32",
   \POS "a32" \ar_{\lambda^{\K\W\widetilde{\R}}(Q, D)} "a33",
   \POS "a33" \ar_{W(Q,\overline{L}d)} "a24",
   \POS "a21" \ar^{LK(\widetilde{R}Q', d)} "a12",
   \POS "a12" \ar^{\lambda^{\K\W\widetilde{\R}}(Q', D')} "a13",
   \POS "a13" \ar^{W(q^{op},\overline{L}D')} "a24"
  \endxy
\end{equation*}

\end{1708091812}

\noindent \emph{Proof}:\\

\noindent First, recall that

\begin{equation*}
K\big((\widetilde{R}q)^{op}, \sim \big) := K_{\widetilde{\R}\Q}\varepsilon^{\J\K\widetilde{R}\Q'}\circ K_{\widetilde{\R}\Q}J_{\widetilde{\R}q}K_{\widetilde{\R}\Q'}\circ \eta^{\K\J\widetilde{R}Q}K_{\widetilde{\R}\Q'}
\end{equation*}\\

\noindent since $J(\mathcal{C}\times\widetilde{R})\dashv_{\hspace{0.4pt}\mathcal{Q}}K(\widetilde{R}^{op}\times\mathcal{D})$. There exists a similar expression for $W(q^{op},\sim)$. Second, the following equation takes place due to the naturality of all of the involved components 

\begin{eqnarray*}
  && W_{\Q}Ld\cdot W_{\Q}\varepsilon^{\V\W\Q'}LD\cdot W_{\Q}V_{q}W_{Q'}LD\cdot \eta^{\W\V\Q}W_{\Q'}LD\cdot \lambda^{\K\W\widetilde{\R}Q}D =\\
 &&\quad\quad\quad\quad W_{\Q}\varepsilon^{\V\W\Q'}LD'\cdot W_{\Q}V_{q}W_{\Q'}LD'\cdot \eta^{\W\V\Q}W_{\Q'}LD'\cdot \lambda^{\K\W\widetilde{\R}Q'}D'\cdot LK_{\widetilde{\R}\Q'}d
\end{eqnarray*}

Therefore, it is left to prove that the following equation takes place

\begin{eqnarray*}
  && \lambda^{\K\W\widetilde{\R}\Q}D\cdot LK_{\widetilde{\R}\Q}\varepsilon^{\J\K\widetilde{\R}\Q'}\cdot LK_{\widetilde{\R}\Q}J_{\widetilde{\R}q}K_{\widetilde{\R}\Q'}D\cdot L\eta^{\K\!\J\widetilde{\R}\Q}K_{\widetilde{\R}\Q'}D\\
  &&\qquad= W_{\Q}\varepsilon^{\V\W\Q'}\overline{L}D\cdot W_{\Q}V_{q}W_{\Q'}\overline{L}D\cdot\eta^{\W\V\Q}W_{\Q'}\overline{L}D\cdot\lambda^{\K\W\widetilde{\R}\Q'}D
\end{eqnarray*}

This is done by the following process:

\begin{eqnarray*}
  &&\lambda^{\K\W\widetilde{\R}\Q}D\cdot LK_{\widetilde{\R}\Q}\varepsilon^{\J\!\K\widetilde{\R}\Q'}\cdot LK_{\widetilde{\R}\Q}J_{\widetilde{\R}q}K_{\widetilde{\R}\Q'}D\cdot L\eta^{\K\!\J\widetilde{\R}\Q}K_{\widetilde{\R}\Q'}D\\
  &&\qquad\quad=\lambda^{\K\W\widetilde{\R}\Q}D\cdot LK_{\widetilde{\R}\Q}\varepsilon^{\J\!\K\widetilde{\R}\Q'}\cdot LK_{(\widetilde{\R}q)^{op}}J_{\widetilde{\R}\Q'}K_{\widetilde{\R}\Q'}D\cdot L\eta^{\K\!\J\widetilde{\R}\Q'}K_{\widetilde{\R}\Q'}D\\
  &&\qquad\quad=\lambda^{\K\W\widetilde{\R}\Q}D\cdot LK_{(\widetilde{\R}q)^{op}}D\cdot K_{\widetilde{\R}\Q'}\varepsilon^{\J\!\K\widetilde{\R}\Q'}D\cdot\eta^{\K\!\J\widetilde{\R}\Q'} K_{\widetilde{\R}\Q'}D\\
  &&\qquad\quad=\lambda^{\K\W\widetilde{\R}\Q}D\cdot LK_{(\widetilde{\R}q)^{op}}D = W_{q^{op}}\overline{L}D\cdot \lambda^{\K\W\widetilde{\R}\Q'}D\\
  &&\qquad\quad=W_{\Q}\varepsilon^{\V\W\Q}\overline{L}D\cdot \eta^{\W\V\Q}W_{\Q}\overline{L}D\cdot W_{q^{op}}\overline{L}D\cdot \lambda^{\K\W\widetilde{\R}\Q'}D\\
  &&\qquad\quad=W_{\Q}\varepsilon^{\V\W\Q}\overline{L}D\cdot W_{\Q}V_{\Q}W_{q^{op}}\overline{L}D\cdot \eta^{\W\V\Q}W_{\Q'}\overline{L}D\cdot \lambda^{\K\W\widetilde{\R}\Q'}D\\
  &&\qquad\quad=W_{\Q}\varepsilon^{\V\W\Q'}\overline{L}D\cdot W_{\Q}V_{q}W_{\Q'}\overline{L}D\cdot\eta^{\W\V\Q}W_{\Q'}\overline{L}D\cdot \lambda^{\K\W\widetilde{\R}\Q'}D 
\end{eqnarray*}

The first equality takes place due to the Proposition \ref{1709030912} where $J_{\widetilde{\R}q}$ and $K_{({\widetilde{\R}q})^{op}}$ are conjugate morphisms. The third one, uses the triangular identity associated to $K_{\widetilde{\R}\Q}$. The fourth one, is due to the fact that $(V_{q}, W_{q^{op}})$ is a 2-cell in $\Adj_{\R}(\mathfrak{C})$. The fifth one uses the triangular identity associated to $W_{\Q}$. The seventh is related to the fact that $V_{q}$ and $W_{q^{op}}$ are conjugate. The rest of the equalities have to do with an involved naturality and therefore the details are spare for those.

\begin{flushright}
$\Box$
\end{flushright}

The mate of $\lambda^{\K\W\!\widetilde{R}\Q}$, $\rho^{\K\W\!\widetilde{R}\Q}$, and the inverse of this last one $\varrho^{\K\W\!\widetilde{R}\Q}$ can also be extended to a dinatural transformation.

\newtheorem{1708122030}[1708081705]{Corollary}
\begin{1708122030}

The transformation defined, for $(Q,Y)$ in $\mathcal{Q}^{op}\times\mathcal{Y}$, as

\begin{equation*}
\varrho^{\K\W\widetilde{R}}(Q,Y) := \varrho^{\K\W\widetilde{R}\Q}Y
\end{equation*}

\noindent is dinatural, \emph{i.e.} for any $(q^{op},y)$ in $\mathcal{Q}^{op}\times\mathcal{Y}$, the following diagram commutes

\begin{equation*}
    \xy<1cm, 1cm>
    \POS (25,0) *+{RW(Q',Y')} = "a12",
    \POS (70,0) *+{K(\widetilde{R}Q', \overline{R}Y')} = "a13",
    \POS (0,-15) *+{RW(Q',Y)} = "a21",
    \POS (95,-15) *+{K(\widetilde{R}Q, \overline{R}Y')} = "a24",
    \POS (25,-30) *+{RW(Q, Y)} = "a32",
    \POS (70,-30) *+{K(\widetilde{R}Q, \overline{R}Y)} = "a33",
    \POS "a21" \ar_{RW(q^{op}, Y)\phantom{h}} "a32",
    \POS "a32" \ar_{\varrho^{\K\W\widetilde{\R}}(Q, Y)} "a33",
    \POS "a33" \ar_{\phantom{h}K(\widetilde{R}Q, \overline{R}y)} "a24",
    \POS "a21" \ar^{RW(Q',\:y)\phantom{h}} "a12",
    \POS "a12" \ar^{\varrho^{\K\W\widetilde{\R}}(Q', Y')} "a13",
    \POS "a13" \ar^{\phantom{h}K((\widetilde{R}q)^{op}, \overline{R}Y')} "a24"
    \endxy
\end{equation*}

\end{1708122030}

\begin{flushright}
$\Box$
\end{flushright}

The inverse of the mate of the dinatural transformation $\lambda^{\K\W\!\widetilde{\R}}$ is the dinatural transformation $\varrho^{\K\W\!\widetilde{\R}}$ then the following is a 1-cell in $\Adj_{\R}(\mathfrak{C})$

\begin{equation*}
(K(\widetilde{R}^{op}\times\mathcal{D}), W, \lambda^{\K\W\!\widetilde{\R}}): \mathcal{Q}^{op}\times\overline{L}\dashv \mathcal{Q}^{op}\times \overline{R}\longrightarrow L\dashv R
\end{equation*}

Therefore, using a left Hopf 1-cell, an object similar to a parametric adjunction could be obtained. This result is summarized, and the corresponding process, into the following statement and definition.\\

\newtheorem{1709121524}[1708081705]{Theorem}

\begin{1709121524}\label{1709121524}

Consider a left Hopf 1-cell of the form

\begin{equation*}
(J, V, \lambda^{\J\V}): L\times\widetilde{L}\dashv R\times\widetilde{R}\longrightarrow \overline{L}\dashv \overline{R}
\end{equation*}

$\phantom{br}$

\noindent and a pair of classical parametric adjunctions $J\dashv_{\hspace{0.4pt}\mathcal{P}} K$ and $V\dashv_{\mathcal{Q}} W$. Then we have\\

\begin{enumerate}

 \item [i)] $(J(\mathcal{C}\times \widetilde{R}), V, \lambda^{\J\V\!\widetilde{\R}}): L\times\mathcal{Q}\dashv R\times\mathcal{Q}\longrightarrow \overline{L}\dashv\overline{R}$.

 \item [ii)] $(K(\widetilde{R}^{op}\times\mathcal{D}), W, \lambda^{\K\W\!\widetilde{\R}}):\mathcal{Q}^{op}\times\overline{L}\dashv \mathcal{Q}^{op}\times \overline{R}\longrightarrow L\dashv R$.

\end{enumerate}

$\phantom{br}$

\noindent as 1-cells in $\Adj_{\R}(\mathfrak{C})$ and for each $Q$ in $\mathcal{Q}$, an adjoint object

 \begin{equation*}
(J_{\widetilde{\R}\Q}, V_{\Q}, \lambda^{\J\V\!\widetilde{\R}\Q})\dashv(K_{\widetilde{\R}\Q}, W_{\Q}, \lambda^{\K\W\!\widetilde{\R}\Q})
 \end{equation*}

 $\phantom{br}$

Therefore, this structure might be defined as a \emph{Hopf parametric adjoint object}, in $\Adj_{\R}(\mathfrak{C})$, and denoted as

\begin{equation*}
\big( J(\mathcal{C}\times\widetilde{R}), V, \lambda^{\J\V\!\widetilde{\R}} \big) \dashv_{\mathcal{Q}} \big( K(\widetilde{R}^{op}\times\mathcal{D}), W, \lambda^{\K\W\!\widetilde{\R}} \big)
\end{equation*}

\end{1709121524}

\begin{flushright}
$\Box$
\end{flushright}


\subsection{The Antipode}


Similar to the definition of an antipode in \cite{bral_homm}, the following corollary is stated in order to define this concept for a left Hopf 1-cell.

\newtheorem{1708131019}{Corollary}[subsection]
\begin{1708131019}

The following transformation

\begin{equation*}
\psi^{\K\!\widetilde{\R}}:= \varrho^{\K\W\!\widetilde{R}}\cdot_{p} \lambda^{\K\W\!\widetilde{R}} = \varrho^{\K\W\!\widetilde{R}}(\mathcal{Q}^{op}\times\overline{L})\circ R\lambda^{\K\W\!\widetilde{R}}:RLK(\widetilde{R}^{op}\times\mathcal{D})\longrightarrow K(\widetilde{R}^{op}\times \overline{R}\overline{L})
\end{equation*}

\noindent is dinatural.

\end{1708131019}

According to \cite{bral_homm}, there is a certain bijection of dinatural transformations, which is now rewritten in this context for a left Hopf 1-cell in $\Adj_{\R}(\mathfrak{C})$.

\newtheorem{1708301657}[1708131019]{Proposition}
\begin{1708301657}\label{1708301657}

There is a bijection between the following dinatural transformations

\begin{enumerate}

\item [i)] $\psi^{\K\!\widetilde{\R}}:RLK(\widetilde{R}^{op}\times\mathcal{D})\longrightarrow K(\widetilde{R}^{op}\times \overline{R}\overline{L}):\mathcal{Q}^{op}\times\mathcal{D}\longrightarrow\mathcal{C}$.

\item [ii)] $\sigma^{\K\!\widetilde{\R}}: RLK(\widetilde{R}^{op}\widetilde{L}^{op}\times\mathcal{D})\longrightarrow K(\mathcal{P}^{op}\times\overline{RL}):\mathcal{P}^{op}\times\mathcal{D}\longrightarrow \mathcal{C}$.

\end{enumerate}

This last dinatural transformation is called \emph{antipode}. 

\begin{flushright}
$\Box$
\end{flushright}

\end{1708301657}


\section{Left Hopf 1-cells and Parametric Adjoint Objects in $\Mnd(\mathfrak{C})$}


In this chapter, the definitions made in the previous section are recalled but this time monads are used. The objective, in this section as in the whole article, is to give an extension of a classical parametric adjunction $J\dashv_{\hspace{0.5pt}\mathcal{P}}K$ within the 2-categorical context of $\Mnd(\mathfrak{C})$.\\

\subsection{Hopf 1-cells}


Consider for this case the following 0-cells $(\mathcal{C}, S)$, $(\mathcal{D}, T)$ and $(\mathcal{P}, E)$. For any functor $J:\mathcal{C}\times\mathcal{P}\longrightarrow \mathcal{D}$ one can think of a 1-cell, in $\Mnd(\mathfrak{C})$, of the form

\begin{equation*}
(J, \psi^{\J}): (\mathcal{C}\times\mathcal{P}, S\times E)\longrightarrow (\mathcal{D}, T)
\end{equation*}

$\phantom{br}$

If one wishes to construct a parametric adjoint object then there must exist a functor $K:\mathcal{P}^{op}\times\mathcal{D}\longrightarrow\mathcal{C}$ that can be extended to a 1-cell in $\Mnd(\mathfrak{C})$, but such an extension presents a problem. Since $E$ is a monad on $\mathcal{P}$, $E^{op}$ is a comonad on $\mathcal{P}^{op}$, therefore  it cannot be proposed a 1-cell $(K, \psi^{\K}):(\mathcal{P}^{op}\times\mathcal{D}, E^{op}\times T)\longrightarrow (\mathcal{C}, S)$, neither in $\Mnd(\mathfrak{C})$ or in the comonad dual of $\Mnd(\mathfrak{C})$ in order to complete a possible parametric adjunction. In the same way as before, a modification of the functors $J$ and $K$ has to be made in order to achive the proposed objective.\\

\newtheorem{1708311910}{Definition}[subsection]
\begin{1708311910}\label{1708311910}

Consider a 1-cell, in $\Mnd(\mathfrak{C})$, of the form

\begin{equation*}
(J,\psi^{\J}):(\mathcal{C}\times\mathcal{P}, S\times E)\longrightarrow (\mathcal{D}, T)
\end{equation*}

\noindent The \emph{left Hopf operator} on $(J,\psi^{\J})$ is the following 1-cell in $\Mnd(\mathfrak{C})$

\begin{equation*}
\mathfrak{H}(J, \psi^{\J}):= (J(\mathcal{C}\times U^{\E}), H(\psi^{\J})) = (J(\mathcal{C}\times U^{\E}), S\times \mathcal{P}^{\E})\longrightarrow (\mathcal{D}, T)
\end{equation*}

\noindent where $H(\psi^{\J})$ is the following natural tranformation

\begin{equation*}
 \xy<1cm,1cm>
  \POS (-5,0) *+{\mathcal{C}\times\mathcal{P}^{\E}} = "a11",
  \POS (25,0) *+{\mathcal{C}\times\mathcal{P}} = "a12",
  \POS (55,0) *+{\mathcal{D}} = "a13",
  \POS (25,-25) *+{\mathcal{C}\times\mathcal{P}} = "a22",
  \POS (55,-25) *+{\mathcal{D}} = "a23",
  \POS (53,-15) *+{} = "x1",
  \POS (40,-23) *+{} = "x2",
  \POS (25,-8) *+{} = "y1",
  \POS (13,-13) *+{} = "y2",
  \POS "a11" \ar_{S\times U^{\E}} "a22",
  \POS "a22" \ar_-{J} "a23",
  \POS "a11" \ar^{\mathcal{C}\times U^{\E}} "a12",
  \POS "a12" \ar^{S\times E} "a22",
  \POS "a12" \ar^-{J} "a13",
  \POS "a13" \ar^{T} "a23",
  \POS "x1" \ar@/_0.5pc/_{\psi^{\J}} "x2",
  \POS "y1" \ar@/_0.4pc/_{S\times U^{\E}\varepsilon^{\F\U\!\E}} "y2",
 \endxy
\end{equation*}

\noindent that is to say,

\begin{equation*}
H(\psi^{\J}) = (S\times U^{\E}\varepsilon^{\F\U\!\E})\cdot_{p}\hspace{0.4pt}\psi^{\J} = J(S\times U^{\E}\varepsilon^{\F\U\!\E})\circ \psi^{\J}(\mathcal{C}\times U^{\E})
\end{equation*}

$\phantom{br}$

\end{1708311910}

\newtheorem{1710241316}[1708311910]{Definition}
\begin{1710241316}

Consider a 1-cell $(J,\psi^{\J}):(\mathcal{C}\times\mathcal{P}, S\times E)\longrightarrow (\mathcal{D}, T)$ in $\Mnd(\mathfrak{C})$. The \emph{left fusion operator} is the following 1-cell in $\Mnd(\mathfrak{C})$

\begin{equation*}
\mathfrak{F}(J, \psi^{\J}):= (J(\mathcal{C}\times E), F(\psi^{\J})) = (J(\mathcal{C}\times E), S\times \mathcal{P})\longrightarrow (\mathcal{D}, T)
\end{equation*}

\noindent where $F(\psi^{\J})$  is the following natural transformation 

\begin{equation*}
 \xy<1cm,1cm>
  \POS (0,0) *+{\mathcal{C}\times\mathcal{P}} = "a11",
  \POS (25,0) *+{\mathcal{C}\times\mathcal{P}} = "a12",
  \POS (55,0) *+{\mathcal{D}} = "a13",
  \POS (25,-25) *+{\mathcal{C}\times\mathcal{P}} = "a22",
  \POS (55,-25) *+{\mathcal{D}} = "a23",
  \POS (53,-15) *+{} = "x1",
  \POS (40,-23) *+{} = "x2",
  \POS (25,-8) *+{} = "y1",
  \POS (15,-15) *+{} = "y2",
  \POS "a11" \ar_{S\times E} "a22",
  \POS "a22" \ar_-{J} "a23",
  \POS "a11" \ar^{\mathcal{C}\times E} "a12",
  \POS "a12" \ar^{S\times E} "a22",
  \POS "a12" \ar^-{J} "a13",
  \POS "a13" \ar^{T} "a23",
  \POS "x1" \ar@/_0.5pc/_{\psi^{\J}} "x2",
  \POS "y1" \ar@/_0.4pc/_{S\times\mu^{\E}} "y2",
 \endxy
\end{equation*}

\noindent that is to say

\begin{equation*}
F(\psi^{\J}) = (S\times \mu^{\E})\cdot_{p}\hspace{0.4pt}\psi^{\J} = J(S\times\mu^{\E})\circ \psi^{\J}(\mathcal{C}\times E)
\end{equation*}

\end{1710241316}

\newtheorem{1710241339}[1708311910]{Definition}
\begin{1710241339}

A \emph{left Hopf} 1-cell, in $\Mnd(\mathfrak{C})$, is of the form $(J, \psi^{\J}): (\mathcal{C}\times\mathcal{P}, S\times E)\longrightarrow (\mathcal{D}, T)$ such that $H(\psi^{\J})$ is invertible. In such case, the inverse is denoted as $N(\psi^{\J})$.

\end{1710241339}

\newtheorem{1709121656}[1708311910]{Definition}
\begin{1709121656}

  A \emph{left fusion} 1-cell, in $\Mnd(\mathfrak{C})$, $(J, \psi^{\J}): (\mathcal{C}\times\mathcal{P}, S\times E)\longrightarrow (\mathcal{D}, T)$ is such that $F(\psi^{\J})$ is invertible. In such case, the inverse is denoted as $G(\psi^{\J})$.

\end{1709121656}

\textbf{Remark:} Later on, it will be checked that the Hopf and fusion 1-cell will be equivalent, which in turn will ease the difference with the Hopf monad definition on \cite{bral_homm}.\\

Consider a classical parametric adjunction $J\dashv_{\hspace{0.5pt}\mathcal{P}} K$ and a left Hopf 1-cell $(J, \psi^{\J}): (\mathcal{C}\times\mathcal{P}, S\times E)\longrightarrow (\mathcal{D}, T)$. Since $H(\psi^{\J})$ is invertible so is the following natural transformation for any $(M, k_{\M})$.

\begin{equation*}
 \xy<1cm,1cm>
  \POS (0,0) *+{\mathcal{C}} = "a11",
  \POS (30,0) *+{\mathcal{C}\times \mathbf{1}} = "a12",
  \POS (60,0) *+{\mathcal{C}\times\mathcal{P}^{\E}} = "a13",
  \POS (90,0) *+{\mathcal{C}\times\mathcal{P}} = "a14",
  \POS (120,0) *+{\mathcal{D}} = "a15",
  \POS (0,-25) *+{\mathcal{C}} = "a21",
  \POS (30,-25) *+{\mathcal{C}\times \mathbf{1}} = "a22",
  \POS (60,-25) *+{\mathcal{C}\times\mathcal{P}^{\E}} = "a23",
  \POS (90,-25) *+{\mathcal{C}\times\mathcal{P}^{\E}} = "a24",
  \POS (88,-15) *+{} = "x1",
  \POS (75,-23) *+{} = "x2",
  \POS (118,-15) *+{} = "y1",
  \POS (105,-23) *+{} = "y2",
  \POS (120,-25) *+{\mathcal{Y}} = "a25",
  \POS "a11" \ar_{S} "a21",
  \POS "a21" \ar_{\rho^{-1}\mathcal{C}} "a22",
  \POS "a22" \ar_{\mathcal{C}\times E_{\M\!k}} "a23",
  \POS "a23" \ar_{\mathcal{C}\times U^{\E}} "a24",
  \POS "a24" \ar_-{J} "a25",
  \POS "a11" \ar^{\rho^{-1}\mathcal{C}} "a12",
  \POS "a12" \ar^{\mathcal{C}\times E_{\M\!k}} "a13",
  \POS "a13" \ar^{\mathcal{C}\times U^{\E}} "a14",
  \POS "a14" \ar^-{J} "a15",
  \POS "a12" \ar^{S\times\mathbf{1}} "a22",
  \POS "a13" \ar_{S\times\mathcal{P}^{\E}} "a23",
  \POS "a14" \ar^{S\times E} "a24",
  \POS "a15" \ar^{T} "a25",
  \POS "x1" \ar@/_0.5pc/_{S\times U^{\E}\varepsilon^{\F\U\!\E}} "x2",
  \POS "y1" \ar@/_0.5pc/_{\psi^{\J}} "y2",
 \endxy
\end{equation*}

The previous invertible natural transformation is denoted as $\psi^{\J\U\!\E}(M, k_{\M})$ or 

\begin{equation*}
\psi^{\J\!\M}: T J_{\M} \longrightarrow J_{\M}S:\mathcal{C}\longrightarrow\mathcal{D},
\end{equation*}

\noindent \textbf{Remark}: If $J\dashv_{\mathcal{P}} K$ then $J(\mathcal{C}\times U^{\E}) \dashv_{\hspace{0.5pt}\mathcal{P}^{\E}} K(U^{{\E}^{op}}\!\times\mathcal{D})$.\\

For any Eilenberg-Moore algebra $(M, k_{\M})$, there is an adjunction $J_{\M}\dashv K_{\M}$ such that the following is a 1-cell in $\Mnd(\mathfrak{C})$

\begin{equation*}
(J_{\M}, \psi^{\J\!\M}): (\mathcal{C}, S)\longrightarrow (\mathcal{D}, T)
\end{equation*}

$\phantom{br}$

\noindent and $\psi^{\J\!\M}$ is invertible. If the Proposition \ref{1612221117} is applied, an adjoint object in $\Mnd(\mathfrak{C})$ is obtained

\begin{equation*}
(J_{\M}, \psi^{\J\!\M})\dashv (K_{\M}, \psi^{\K\!\M})
\end{equation*}

\noindent where $\psi^{\K\!\M} = ad(\psi^{\J\!\M})$. The last natural transformation can be further extended.

\newtheorem{1709041732}[1708311910]{Proposition}
\begin{1709041732}

The transformation $\psi^{\K\!\M}$, on $(M, k_{\M})$, can be extended to the following dinatural transformation

\begin{equation*}
\psi^{\K\U\!\E}: SK(U^{\E^{op}}\times\mathcal{D})\longrightarrow K(U^{\E^{op}}\times\mathcal{D})((\mathcal{P}^{\E})^{op}\times T): (\mathcal{P}^{\E})^{op}\times\mathcal{D}\longrightarrow \mathcal{C}.
\end{equation*}

\end{1709041732}

$\phantom{br}$

\noindent \emph{Proof:}\\

\noindent Define $\psi^{\K\U\!\E}$ for  $\big((M, k_{\M}), D\big)$, in $(\mathcal{P}^{\E})^{op}\times\mathcal{D}$, as

\begin{equation*}
\psi^{\K\U\!\E}\big((M, k_{\M}), D \big) := \psi^{\K\!\M}D.
\end{equation*}

$\phantom{br}$

\noindent The proof of the commutativity for the corresponding morphism $\overline{p}^{\hspace{0.5pt}op}:(M', k_{\M'})\longrightarrow (M, k_{\M})$ is left to the reader.

\begin{flushright}
$\Box$
\end{flushright}

In the context of the previous Proposition, the dinatural transformation can be denoted as $\psi^{\K\U\!\E} :=  ad(H(\psi^{\J}))$ or $\psi^{\K\U\!\E} : = H^{\sharp}(\psi^{\J})$.\\

The previous proccess can be summarized into the following Theorem.

\newtheorem{1709170010}[1708311910]{Theorem}
\begin{1709170010}\label{1709170010}
  
Consider a classical parametric adjunction $J\dashv_{\hspace{0.5pt}\mathcal{P}}K$ and a left Hopf 1-cell in $\Mnd(\mathfrak{C})$ of the form

\begin{equation*}
(J, \psi^{\J}):(\mathcal{C}\times\mathcal{P}, S\times E)\longrightarrow(\mathcal{D},T)
\end{equation*}

Then the following

\begin{enumerate}

\item [i)] $(J(\mathcal{C}\times U^{\E}), \psi^{\J\U\!\E}): (\mathcal{C}\times\mathcal{P}^{\E}, S\times\mathcal{P}^{\E})\longrightarrow(\mathcal{D}, T)$,

\item [ii)] $(K(U^{\E^{op}}\times\mathcal{D}), \psi^{\K\U\!\E}):\big( (\mathcal{P}^{\E})^{op}\times\mathcal{D}, (\mathcal{P}^{\E})^{op}\times T\big)\longrightarrow (\mathcal{C}, S)$

\end{enumerate}

\noindent are 1-cells in $\Mnd(\mathfrak{C})$ and for each $(M, k_{\M})$ in $\mathcal{P}^{\E}$ there is an adjoint object

\begin{equation*}
(J_{\M}, \psi^{\J\!\M})\dashv (K_{\M}, \psi^{\K\!\M})
\end{equation*}

Therefore, this structure might be defined as a \emph{Hopf parametric adjoint object}, in $\Mnd(\mathfrak{C})$, and denoted as

\begin{equation*}
(J(\mathcal{C}\times U^{\E}), \psi^{\J\U\!\E}) \dashv_{\hspace{0.4pt}\mathcal{P}^{\E}} (K(U^{\E^{op}}\times\mathcal{D}), \psi^{\K\U\!\E})
\end{equation*}

\begin{flushright}
$\Box$
\end{flushright}

\end{1709170010}


\subsection{Antipode}


Analogous to the bijection in Proposition \ref{1708301657},  consider the dinatural transformation

\begin{equation*}
\psi^{\K\U\!\E}:SK(U^{\E^{op}}\times 1_{\mathcal{D}})\longrightarrow K(U^{\E^{op}}\times 1_{\mathcal{D}})(\mathcal{P}^{op}\times T)
\end{equation*}

$\phantom{br}$

\noindent whisker it with the functor $F^{\E^{op}}\!\!\times \mathcal{D}$ and compose it with the natural transformation $K(\varepsilon^{\U\!\F\!\E^{op}}\!\times T)$ to get


\begin{equation*}
\sigma := K(\varepsilon^{\U\!\F\!\E^{op}}\!\times T)\circ \psi^{\K\U\!\E}(F^{\E^{op}}\!\!\times \mathcal{D}): SK(E^{op}\times\mathcal{D})\longrightarrow K(\mathcal{P}^{op}\!\times T): \mathcal{P}^{op}\times \mathcal{D}\longrightarrow\mathcal{C}
\end{equation*}

$\phantom{br}$

\noindent whose component at $(P, D)$ in $\mathcal{P}^{op}\times\mathcal{D}$, noting that $\varepsilon^{\U\!\F\!\E^{op}}P = (\eta^{\F\U\!\E}P)^{op} = (\eta^{\E}P)^{op}$, is

\begin{equation*}
\sigma^{\K}(P, D) = K\big((\eta^{\E}P)^{op}, TD\big)\cdot \psi^{\K\U\!\E}((EP,(\mu^{\E}P)^{op} ), D)
\end{equation*}

$\phantom{br}$

\noindent therefore, $\sigma^{\K}(P, D): SK(EP, D)\longrightarrow K(P, TD)$.

$\phantom{br}$

In \cite{bral_homm}, A. Brugieres et. al. called this natural transformation (left) \emph{antipode}. As pointed out by them, there is a bijection between the dinatural transformations $\psi^{\K\U\!\E}$ and $\sigma^{\K}$, where the inverse of the bijection acts on $\sigma^{\K}$ as follows

\begin{equation*}
\iota := \sigma^{\K}(U^{\E^{op}}\times \mathcal{D})\circ SK(\eta^{\F\U\!\E^{op}}\times\mathcal{D}): SK(U^{\E^{op}}\times\mathcal{D})\longrightarrow K(U^{E^{op}}\times \mathcal{D})
\end{equation*}

$\phantom{br}$

\noindent whose component at the object $\big((M, k_{\M}), D \big)$, noting that $\eta^{\F\U\!\E^{op}}(M,k_{\M}) = (\varepsilon^{\F\U\!\E}(M, k_{\M}))^{op} = (k_{\M})^{op}$, is

\begin{equation*}
\iota\big((M, k_{\M}), D \big) = \sigma^{\K}(M, D)\cdot SK(k^{op}_{\M} , D): SK(M,D)\longrightarrow K(M, TD)
\end{equation*}

$\phantom{br}$

\noindent reminiscent of the properties for $\psi^{\K\U\!\E}$, as a 1-cell in $\Mnd(\mathfrak{C})$, the equations that fulfills this antipode are the following

\begin{eqnarray*}
\sigma^{\K}\circ\mu^{\sS}K(E^{op}\times\mathcal{D}) &=& K(1_{\ast}\times\mu^{\T})\circ \sigma^{\K}(1_{\ast}\times T)\circ S\sigma^{\K}(E^{op}\times\mathcal{D})\circ SSK\big((\mu^{\E})^{op}\mathcal{D}\big)  \\
\sigma^{\K}\circ\eta^{\sS}K(E^{op}\times\mathcal{D}) &=& K(\varepsilon^{\U\F\E^{op}}\times\eta^{\T}) = K\big((\eta^{\E})^{op}\times\eta^{\T} \big) 
\end{eqnarray*}

Compare these equations with those equivalent as in Proposition 3.8.b, \cite{bral_homm}.

\section{Left Hopf 1-cells through the 2-adjunction $\Phi_{\E}^{\sfM}\dashv \Psi_{\E}^{\sfM}$}

\subsection{Comparing Hopf 1-cells}


Consider the 1-cell $(J, V, \lambda^{\J\V}): L\times\widetilde{L}\dashv R\times\widetilde{R}\longrightarrow \overline{L}\dashv\overline{R}$ in $\Adj_{\R}(\mathfrak{C})$. This induces a 1-cell in $\Mnd(\mathfrak{C})$ of the form $\Phi_{\E}^{\sfM}(J, V, \lambda^{\J\V}) = (J, \Phi(\lambda^{\J\V})) = (J, \varrho^{\J\V}(L\times\widetilde{L})\circ \overline{R}\lambda^{\J\V})$, where $\varrho^{\J\V}$ is the inverse of the mate $\rho^{\J\V} = {}_{\R\times\widetilde{R}}\mathbf{m}_{\overline{\R}}(\lambda^{\J\V})$. Therefore

\begin{eqnarray*}
H(\Phi(\lambda^{\J\V}))
  &=& H(\varrho^{\J\V}(L\times\widetilde{L})\circ \overline{R}\lambda^{\J\V})\\
  &=& J(RL\times\widetilde{R}\varepsilon^{\widetilde{\sL}\widetilde{\R}})\circ \varrho^{\J\V}(L\times\widetilde{L})(\mathcal{C}\times\widetilde{R})\circ \overline{R}\lambda^{\J\V}(\mathcal{C}\times\widetilde{R})\\
  &=& \varrho^{\J\V}(L\times\mathcal{Q})\circ\overline{R}\big( V(L\times\varepsilon^{\widetilde{\sL}\widetilde{\R}})\circ \lambda^{\J\V}(\mathcal{C}\times\widetilde{R})  \big)\\
  &=& \varrho^{\J\V}(L\times\mathcal{Q})\circ\overline{R}H(\lambda^{\J\V}) = \Phi(H(\lambda^{\J\V}))
\end{eqnarray*}

\noindent where the last equality takes place since ${}_{\R\times\widetilde{R}}\mathbf{m}_{\overline{\R}}(\lambda^{\J\V}) = {}_{\R\times\mathcal{Q}}\mathbf{m}_{\overline{\R}}(H(\lambda^{\J\V}))$. Then, the following proposition can be stated.

\newtheorem{1711011610}{Proposition}[subsection]
\begin{1711011610}

Consider the 1-cell $(J, V, \lambda^{\J\V}): L\times\widetilde{L}\dashv R\times\widetilde{R}\longrightarrow \overline{L}\dashv\overline{R}$ in $\Adj_{\R}(\mathfrak{C})$, such that $\overline{R}$ reflects isomorphisms, then the following statements are equivalent:

\begin{enumerate}

\item [i)] $H(\lambda^{\J\V})$ is invertible, \emph{i.e.} the 1-cell is left Hopf in $\Adj_{\R}(\mathfrak{C})$.

\item [ii)] $H(\Phi(\lambda^{\J\V}))$ is invertible, \emph{i.e.} the induced 1-cell $\Phi(\lambda^{\J\V})$ is left Hopf in $\Mnd(\mathfrak{C})$.

\end{enumerate}

\end{1711011610}

\noindent \emph{Proof}: \\

\noindent $i)\Rightarrow ii)$ \\

If $H(\lambda^{\J\V})$ is invertible, so is $\varrho^{\J\V}(L\times\mathcal{Q})\circ\overline{R}H(\lambda^{\J\V}) = \Phi(H(\lambda^{\J\V}))$ and the conclusion follows from the previous equality.\\

\noindent $ii)\Rightarrow i)$\\

If $H(\Phi(\lambda^{\J\V}))$ is invertible so is $\rho^{\J\V}(L\times\mathcal{Q})\circ H(\Phi(\lambda^{\J\V})) = \overline{R}H(\lambda^{\J\V})$, since $\overline{R}$ reflects isomorphisms $H(\lambda^{\J\V})$ is invertible.

\begin{flushright}
$\Box$
\end{flushright}

The inverses are related as follows

\begin{equation*}
N(\Phi(\lambda^{\J\V})) = \overline{R}N(\lambda^{\J\V})\circ \rho^{\J\V}(L\times\mathcal{Q})
\end{equation*}

\subsection{Hopf Parametric Adjunctions through the 2-adjunction $\Phi_{\E}^{\sfM}\dashv \Psi_{\E}^{\sfM}$}


Using the unit of the 2-adjunction $\Phi_{\E}^{\sfM}\dashv \Psi_{\E}^{\sfM}$ for the 1-cell $(J(\mathcal{C}\times\widetilde{R}), V, H(\lambda^{\J\V}))$ the following proposition can be stated.

\newtheorem{1801061404}{Proposition}[subsection]
\begin{1801061404}\label{1801061404}

Consider the following list:

\begin{itemize}

\item [i)] $L\dashv R:\mathcal{C}\longrightarrow\mathcal{X}$ and the induced monad $(\mathcal{C}, S)$.

\item [ii)] $\overline{L}\dashv\overline{R}:\mathcal{D}\longrightarrow\mathcal{Y}$ and the induced monad $(\mathcal{D}, T)$.

\item [iii)] $\widetilde{L}\dashv\widetilde{R}:\mathcal{P}\longrightarrow\mathcal{Q}$ and the induced monad $(\mathcal{P}, E)$.
  
\end{itemize}

\noindent and suppose that $(J, V, \lambda^{\J\V})$, is a Hopf 1-cell. Therefore, there exists the following pair of commuting diagramms in $\Adj_{\R}(\mathfrak{C})$

\begin{equation*}
\begin{array}{ccc}
\xy<1cm,1cm>
\POS (0,0) *+{\mathcal{C}\times\mathcal{Q}} = "a11",
\POS (35,0) *+{\mathcal{D}} = "a12",
\POS (0,-25) *+{\mathcal{X}\times\mathcal{Q}} = "a21",
\POS (35,-25) *+{\mathcal{Y}} = "a22",
\POS (23,-19) *+{\mathcal{C}\times\mathcal{Q}} = "b11",
\POS (58,-19) *+{\mathcal{D}} = "b12",
\POS (23,-44) *+{\mathcal{C}^{\sS}\times\mathcal{Q}^{1_{\mathcal{Q}}}} = "b21",
\POS (58,-44) *+{\mathcal{D}^{\T}} = "b22",
\POS "a11" \ar^{J(\mathcal{C}\times\widetilde{R})} "a12",
\POS "a12" \ar@<-2pt>_{\overline{L}} "a22",
\POS "a22" \ar@<-2pt>_{\overline{R}} "a12",
\POS "a11" \ar@<-2pt>_{L\times\mathcal{Q}} "a21",
\POS "a21" \ar@<-2pt>_{R\times\mathcal{Q}} "a11",
\POS "a21" \ar_{V} "a22",
\POS "b11" \ar^{\hspace{0.3cm}J(\mathcal{C}\times\widetilde{R})} "b12",
\POS "b12" \ar@<-2pt>_{F^{\T}} "b22",
\POS "b22" \ar@<-2pt>_{U^{\T}} "b12",
\POS "b11" \ar@<-2pt>_{F^{\sS}\times F^{1_{\mathcal{Q}}}} "b21",
\POS "b21" \ar@<-2pt>_{U^{\sS}\times U^{1_{\mathcal{Q}}}} "b11",
\POS "b21" \ar_{[J(\mathcal{C}\times\widetilde{R})]^{H\psi}} "b22",
\POS "a11" \ar "b11",
\POS "a12" \ar "b12",
\POS "a21" \ar_{\mathfrak{K}_{\mathcal{X}}^{\sS}\times\mathfrak{K}_{\mathcal{Q}}^{1_{\mathcal{Q}}}} "b21",
\POS "a22" \ar_{\mathfrak{K}_{\mathcal{Y}}^{\T}} "b22"
\endxy &&
\xy<1cm,1cm>
\POS (0,0) *+{\mathcal{Q}^{op}\times\mathcal{D}} = "a11",
\POS (35,0) *+{\mathcal{C}} = "a12",
\POS (0,-25) *+{\mathcal{Q}^{op}\times\mathcal{Y}} = "a21",
\POS (35,-25) *+{\mathcal{X}} = "a22",
\POS (23,-19) *+{\mathcal{Q}^{op}\times\mathcal{D}} = "b11",
\POS (58,-19) *+{\mathcal{C}} = "b12",
\POS (23,-44) *+{(\mathcal{Q}^{op})^{1_{\mathcal{Q}^{op}}}\times\mathcal{D}^{\T}} = "b21",
\POS (58,-44) *+{\mathcal{C}^{\sS}} = "b22",
\POS "a11" \ar^{K(\widetilde{R}^{op}\times\mathcal{D})} "a12",
\POS "a12" \ar@<-2pt>_{L} "a22",
\POS "a22" \ar@<-2pt>_{R} "a12",
\POS "a11" \ar@<-2pt>_{\mathcal{Q}^{op}\times\overline{L}} "a21",
\POS "a21" \ar@<-2pt>_{\mathcal{Q}^{op}\times\overline{R}} "a11",
\POS "a21" \ar_{W} "a22",
\POS "b11" \ar^{\hspace{0.5cm}K(\widetilde{R}^{op}\times\mathcal{D})} "b12",
\POS "b12" \ar@<-2pt>_{F^{\sS}} "b22",
\POS "b22" \ar@<-2pt>_{U^{\sS}} "b12",
\POS "b11" \ar@<-2pt>_{F^{1_{\ast}}\times F^{\T}} "b21",
\POS "b21" \ar@<-2pt>_{U^{1_{\ast}}\times U^{\T}} "b11",
\POS "b21" \ar_-{[K(\widetilde{R}^{op}\times\mathcal{D})]^{H\!\sharp\psi}} "b22",
\POS "a11" \ar "b11",
\POS "a12" \ar "b12",
\POS "a21" \ar_{\mathfrak{K}_{\mathcal{Y}}^{\T}\times\mathfrak{K}_{\mathcal{Q}^{op}}^{1_{\ast}}} "b21",
\POS "a22" \ar_{\mathfrak{K}_{\mathcal{X}}^{\sS}} "b22"
\endxy 
\end{array}
\end{equation*}


\end{1801061404}

\begin{flushright}
$\Box$
\end{flushright}

Consider the 2-adjunction $\Phi_{\E}^{\sfM}\dashv \Psi_{\E}^{\sfM}$, this structure gives particular classes of isomorphisms of categories. Certain 0-cells are chosen in order to get an adequate isomorphism, for example, consider the following monads, 0-cells in $\Mnd(\mathfrak{C})$, $(\mathcal{C}, S)$, $(\mathcal{P}, E)$ and $(\mathcal{D}, T)$ and construct the 0-cells $F^{\sS}\times F^{\E}\dashv U^{\sS}\times U^{\E}$ and $F^{\T}\dashv U^{\T}$ , in $\Adj_{\R}(\mathfrak{C})$, therefore exists the following isomorphism of categories

\begin{equation*}
Hom_{\Adj_{\R}(\mathfrak{C})}\big(F^{\sS}\times \mathcal{P}^{\E}\dashv U^{\sS}\times \mathcal{P}^{\E}, F^{\T}\dashv U^{\T}\big)
\cong
Hom_{\Mnd(\mathfrak{C})}\big( (\mathcal{C}\times\mathcal{P}^{\E}, S\times \mathcal{P}^{\E}), (\mathcal{D}, T) \big)\\
\end{equation*}

$\phantom{br}$

\noindent Similar isomorphisms exists for combinations of the 0-cells $(\mathcal{P}^{\E})^{op}\times F^{\T}\dashv (\mathcal{P}^{\E})^{op}\times U^{\T}$, $F^{\T}\dashv U^{\T}$ and $F^{\sS}\dashv U^{\sS}$.\\

The following theorem can come forth which combines, through the 2-adjunction $\Phi_{\E}^{\sfM}\dashv \Psi_{\E}^{\sfM}$ and the corresponding isomorphisms, the parametric adjoint objects in Theorem \ref{1709121524} and Theorem \ref{1709170010}.

\newtheorem{1711011826}[1801061404]{Theorem}
\begin{1711011826}\label{1711011826}
 Consider a left Hopf 1-cell $(J, \psi^{\J}): (\mathcal{C}\times\mathcal{P}, S\times E)\longrightarrow (\mathcal{D}, T)$ in $\Mnd(\mathfrak{C})$ whose functor is part of a classical parametric adjunction $J\dashv_{\hspace{0.5pt}\mathcal{P}}K$. Therefore, there exists a bijection between the following structures

\begin{enumerate}

\item [i)] Hopf parametric adjunctions, in $\Adj_{\R}(\mathfrak{C})$, of the form

\begin{equation*}
\Big( J(\mathcal{C}\times U^{\E}) , \big[J(\mathcal{C}\times U^{\E})\big]^{\sH\psi}, \lambda^{\J\!\sH\psi}\Big)\dashv_{\hspace{0.5pt}\mathcal{P}^{\E}} \Big( K(U^{\E^{op}}\times\mathcal{D}) , \big[K(U^{\E^{op}}\times\mathcal{D})\big]^{\sH \sharp\psi}, \lambda^{\K\!\sH\sharp\psi}\Big)
\end{equation*}

\item [ii)] Hopf parametric adjunctions, in $\Mnd(\mathfrak{C})$, of the form

\begin{equation*}
\big( J(\mathcal{C}\times U^{\E}), H(\psi^{\J})\big)\dashv_{\hspace{0.5pt}\mathcal{P}^{\E}}\big( K(U^{\E^{op}}\times\mathcal{D}), H^{\sharp}(\psi^{\J})\big)
\end{equation*}

\end{enumerate}

\end{1711011826}


\section{Lifting parametric adjunctions}


In order to lift a classical parametric adjunction, to some Eilenberg-Moore categories of algebras, there are some further discussion and calculations to be done.

\subsection{Hopf and Fusion 1-cells}


Consider the 1-cell $(J, V, \lambda^{\J\V}):L\times\widetilde{L}\dashv R\times \widetilde{R}\longrightarrow \overline{L}\dashv\overline{R}$. Similar to the relation of the Hopf operators, there is the following relation between the fusion and Hopf operators

\begin{equation}\label{1801051245}
F(\Phi(\lambda^{\J\V})) = \big(\varrho^{\J\V}\cdot_{p} H(\lambda^{\J\V})\big)(\mathcal{C}\times\widetilde{L})
\end{equation}

The following lemma is required.

\newtheorem{1801051438}{Lemma}[subsection]
\begin{1801051438}
Given an adjunction of the form $\widetilde{L}\dashv\widetilde{R}:\mathcal{P}\longrightarrow\mathcal{Q}$ and a natural transformation $\alpha: A\widetilde{R}\longrightarrow B\widetilde{R}$, where $A$ and $B$ are arbitrary parallel functors, with domain $\mathcal{Q}$. Therefore, $\alpha$ is invertible if $\alpha\widetilde{L}$ is so. In this case the inverse of the component $\alpha Q$ is the following

\begin{equation}
\alpha^{-1} Q = A\widetilde{R}\varepsilon^{\widetilde{\sL}\widetilde{\R}}Q\cdot (\alpha\widetilde{L})^{-1}\widetilde{R}Q\cdot B\eta^{\widetilde{\R}\widetilde{\sL}}\widetilde{R}Q 
\end{equation}

\end{1801051438}

\noindent \emph{Proof}:\\

The proof is similar to the Lemma 2.19 given in \cite{bral_homm} but this time the following split fork is used

\begin{equation*}
\xy<1cm,1cm>
\POS (0,0) *+{\widetilde{R}\widetilde{L}\widetilde{R}\widetilde{L}\widetilde{R}Q} = "a11",
\POS (35,0) *+{\widetilde{R}\widetilde{L}\widetilde{R}Q} = "a12",
\POS (60,0) *+{\widetilde{R}Q} = "a13",
\POS "a11" \ar@<-2pt>_{\widetilde{R}\widetilde{L}\widetilde{R}\varepsilon^{\widetilde{\sL}\!\widetilde{\R}}Q} "a12",
\POS "a11" \ar@<2pt>^{\widetilde{R}\varepsilon^{\widetilde{\sL}\!\widetilde{\R}}\widetilde{L}\widetilde{R}Q} "a12",
\POS "a12" \ar^{\widetilde{R}\varepsilon^{\widetilde{\sL}\!\widetilde{\R}}Q} "a13",
\endxy
\end{equation*}

\begin{flushright}
$\Box$
\end{flushright}

The following proposition can be written.

\newtheorem{1711050055}[1801051438]{Proposition}
\begin{1711050055}
  
Consider the 1-cell $(J, V, \lambda^{\J\V}):L\times\widetilde{L}\dashv R\times \widetilde{R}\longrightarrow \overline{L}\dashv\overline{R}$, such that $\overline{R}$ reflects isomorphisms, then the following statements are equivalent

\begin{enumerate}
  
\item [i)] $H(\lambda^{\J\V})$ is invertible, \emph{i.e.} the 1-cell is left Hopf in $\Adj_{\R}(\mathfrak{C})$.

\item [ii)] $F(\Phi(\lambda^{\J\V}))$ is invertible, \emph{i.e.} the 1-cell is left fusion in $\Mnd(\mathfrak{C})$.
  
\end{enumerate}

\end{1711050055}

\noindent \emph{Proof}:\\

\noindent $i) \Rightarrow ii)$ Clear by taking into consideration \eqref{1801051245}.\\

\noindent $ii) \Rightarrow i)$ In order to use the previous lemma, take $C$ in $\mathcal{C}$ and the natural transformation $\alpha$ is given by

\begin{equation*}
\big(\varrho^{\J\V}\cdot_{p} H(\lambda^{\J\V})\big)(C,\sim): \overline{R}\overline{L} J(C,\widetilde{R})\longrightarrow J(RLC,\widetilde{R}),
\end{equation*}

\noindent therefore $\varrho^{\J\V}\cdot_{p} H(\lambda^{\J\V})$ is invertible and so is $\overline{R}H(\lambda^{\J\V})$, since $\overline{R}$ reflects isomorphisms $H(\lambda^{\J\V})$ is also invertible.

\begin{flushright}
$\Box$
\end{flushright}

The inverses are related as follows

\begin{equation*}
\overline{R}N(\lambda^{\J\V})\circ \rho^{\J\V} (L\times\mathcal{Q}) =  \overline{R}\overline{L}J(\mathcal{C}\times \widetilde{R}\varepsilon^{\widetilde{\sL}\!\widetilde{\R}})\circ G(\Phi(\lambda^{\J\V}))(\mathcal{C}\times \widetilde{R}) \circ J(RL\times\eta^{\widetilde{\R}\!\widetilde{\sL}}\widetilde{R})
\end{equation*}

The reader is compelled to check the same expression in Lemma 2.18 \cite{bral_homm}.

\newtheorem{1711050101}[1801051438]{Corollary}
\begin{1711050101}

Consider a 1-cell $(J, \psi^{\J}):(\mathcal{C}\times\mathcal{P}, S\times E)\longrightarrow (\mathcal{D}, T)$. Therefore $F(\psi^{\J})$ is invertible iff $H(\Psi(\psi^{\J}))$ is invertible, \emph{i.e.} the 1-cell is Hopf iff is fusionable.

\end{1711050101}

This corollary allows the author to keep using the adjective Hopf without losing the generality of the results.


\subsection{Hopf Parametric Liftings}


Without any further ado, the main Theorem of the article is stated and proved. \\

\newtheorem{1711061222}{Theorem}[subsection]
\begin{1711061222}

Consider a parametric adjunction $J\dashv_{\hspace{0.5pt}\mathcal{P}} K$, and 0-cells in $\Mnd(\mathfrak{C})$ of the form $(\mathcal{C}, S)$, $(\mathcal{D}, T)$ and $(\mathcal{P}, E)$. There is a bijection between the following structures

\begin{enumerate}

\item [i)] Parametric Adjoint Liftings $\widehat{J}\dashv_{\hspace{0.5pt}\mathcal{P}^{\E}} \widehat{K}$, where $(J, \widehat{J}, \lambda^{\J\!\widehat{J}})$ is a Hopf 1-cell in $\Adj_{\R}(\mathfrak{C})$. This lifted parametric adjunction makes the following diagrams commutative

\begin{equation*}
\begin{array}{ccc}
\xy<1cm,1cm>
  \POS (0,0) *+{\mathcal{C}\times\mathcal{P}} = "a11",
  \POS (25,0) *+{\mathcal{D}} = "a12",
  \POS (0,-20) *+{\mathcal{C}^{\sS}\times\mathcal{P}^{\E}} = "a21",
  \POS (25,-20) *+{\mathcal{D}^{\T}} = "a22",
  \POS "a21" \ar^{U^{\sS}\times U^{\E}} "a11",
  \POS "a11" \ar^-{J} "a12",
  \POS "a21" \ar_-{\widehat{J}} "a22",
  \POS "a22" \ar_{U^{\T}} "a12",
  \endxy & &
  \xy<1cm,1cm>
  \POS (0,0) *+{\mathcal{P}^{op}\times\mathcal{D}} = "a11",
  \POS (25,0) *+{\mathcal{C}} = "a12",
  \POS (0,-20) *+{(\mathcal{P}^{\E})^{op}\times\mathcal{D}^{\T}} = "a21",
  \POS (25,-20) *+{\mathcal{C}^{\sS}} = "a22",
  \POS "a21" \ar^{U^{\E^{op}}\times U^{\T}} "a11",
  \POS "a11" \ar^-{K} "a12",
  \POS "a21" \ar_-{\widehat{K}} "a22",
  \POS "a22" \ar_{U^{\sS}} "a12",
\endxy 
\end{array}
\end{equation*}

\item [ii)] Hopf parametric adjunctions of the form

\begin{equation*}
(J(\mathcal{C}\times U^{\E}), H(\psi^{\J})) \dashv_{\hspace{0.5pt}\mathcal{P}^{\E}} \big( K(U^{\E^{op}}\times\mathcal{D}), \psi^{\K\U\!\E}\big)
\end{equation*}

\noindent where $\psi^{\K\U\!\E} = ad(H(\psi^{\J}))$

\end{enumerate}
\end{1711061222}

\noindent \emph{Proof}: \\



Induce the following left Hopf 1-cell $(J, \psi^{\J}):(\mathcal{C}\times\mathcal{P}, S\times E)\longrightarrow(\mathcal{D}, T)$ in $\Mnd(\mathcal{C})$, where 
$\psi^{\J} := \Phi(\lambda^{\J\widehat{\J}})$. According to Theorem \ref{1711011826} there is a bijection between 

\begin{enumerate}

\item [$\cdot$)] Hopf parametric adjunctions in $\Mnd(\mathfrak{C})$ of the form

\begin{equation*}
\big( J(\mathcal{C}\times U^{\E}), H(\psi^{\J})\big)\dashv_{\hspace{0.5pt}\mathcal{P}^{\E}}\big( K( U^{\E^{op}} \times \mathcal{D}), H^{\sharp}(\psi^{\J})\big)
\end{equation*}

\item [$\cdot$)] Hopf parametric adjunctions in $\Adj_{\R}(\mathfrak{C})$ of the form

\begin{equation*}
\big( J(\mathcal{C}\times U^{\E}) , [J(\mathcal{C}\times U^{\E})]^{\sH\psi}, \lambda^{\J\!\sH\psi}\big)\dashv_{\hspace{0.5pt}\mathcal{P}^{\E}} \big( K(U^{\E^{op}}\times\mathcal{D}) , [K(U^{\E^{op}}\times\mathcal{D})]^{\sH \sharp\psi}, \lambda^{\K\!\sH\sharp\psi}\big)
\end{equation*}

\end{enumerate}

$\phantom{br}$

Taking into account this result, the bijection can be given as follows.\\

The first lifting diagram can be seen as a 1-cell $(J, \widehat{J}, \lambda^{\J\widehat{\J}})$ in $\Adj_{\R}(\mathfrak{C})$ and the following 1-cell can be constructed  $(J(\mathcal{C}\times U^{\E}), \widehat{J}, H(\lambda^{\J\widehat{J}})):F^{\sS}\times\mathcal{P}^{\E}\dashv U^{\sS}\times \mathcal{P}^{\E}\longrightarrow F^{\T}\dashv U^{\T}$. Using the Proposition \ref{1801061404} with this last 1-cell the commutative diagram can be obtained

\begin{equation*}
 \xy<1cm,1cm>
  \POS (0, 0) *+{\mathcal{C}^{\sS}\times\mathcal{P}^{\E}} = "a11",
  \POS (40,0) *+{\mathcal{D}^{\T}} = "a12",
  \POS (0,-25) *+{\mathcal{C}^{\sS}\times (\mathcal{P}^{\E})^{1_{\ast}}} = "a21",
  \POS (40,-25) *+{\mathcal{D}^{\T}} = "a22",
  \POS "a11" \ar_{1_{\mathcal{C}^{\sS}}\times F^{1_{\ast}}} "a21",
  \POS "a21" \ar_-{[J(\mathcal{C}\times U^{\E})]^{H\psi}} "a22",
  \POS "a11" \ar^{\widehat{J}} "a12",
  \POS "a12" \ar^{1_{\mathcal{D}^{\T}}} "a22",
 \endxy
\end{equation*}

A similar argument, for the 1-cell $(K(U^{{\E}^{op}}\times \mathcal{D}), \widehat{K})$, gives the following commutative diagram

\begin{equation*}
  \xy<1cm,1cm>
  \POS (0, 0) *+{(\mathcal{P}^{\E})^{op}\times \mathcal{D}^{\T}} = "a11",
  \POS (40,0) *+{\mathcal{C}^{\sS}} = "a12",
  \POS (0,-25) *+{(\mathcal{P}^{\E^{op}})^{1_{\ast}}\times\mathcal{C}^{\sS}} = "a21",
  \POS (40,-25) *+{\mathcal{C}^{\sS}} = "a22",
  \POS "a11" \ar_{F^{1_{\ast}}\times 1_{\mathcal{D}^{\T}}} "a21",
  \POS "a21" \ar_-{[K(U^{\E^{op}}\times\mathcal{D})]^{H\sharp\psi}} "a22",
  \POS "a11" \ar^{\widehat{K}} "a12",
  \POS "a12" \ar^{1_{\mathcal{C}^{\sS}}} "a22",
  \endxy
\end{equation*}

These two last diagrams give the bijection of the forgetful diagrams with the components of the Hopf parametric adjunction in $\Adj_{\R}(\mathfrak{C})$.

\begin{flushright}
$\Box$
\end{flushright}

\subsection{Hopf monads on monoidal categories}

Consider the case a closed monoidal category $(\mathcal{C}, \otimes, \Box, I)$. In \cite{lojt_apke}, there is a bijection between monoidal liftings, 1-cells in $(\otimes, \widehat{\otimes}, \lambda^{\otimes\widehat{\otimes}}): (F^{\sS}\times F^{\sS}, U^{\sS}\times U^{\sS})\longrightarrow (F^{\sS}\dashv U^{\sS})$  in $\Adj_{\R}(\mathfrak{C})$, and opmonoidal monads, 1-cells $(\otimes, \psi^{\otimes}):(\mathcal{C}\times\mathcal{C}, S\times S)\longrightarrow S$ in $\Mnd(\mathfrak{C})$.\\

If the closure functor is to be lifted when $(\otimes, \widehat{\otimes}, \lambda^{\otimes\widehat{\otimes}})$ is a Hopf 1-cell, the previous calculations show that 

\begin{eqnarray*}
&& (J(\mathcal{C}\times U^{\E}), [J(\mathcal{C}\times U^{\E})]^{H\psi}, \lambda^{\J\!\sH\psi}) = (\otimes(\mathcal{C}\times U^{\sS}), [\otimes(\mathcal{C}\times U^{\sS})]^{\sH\psi}, \lambda^{\otimes\sH\psi})\cong (\otimes, \widehat{\otimes}, \lambda^{\otimes\widehat{\otimes}})\\
&& (K(U^{\E^{op}}\times\mathcal{D}), [K(U^{\E^{op}}\times\mathcal{D})]^{H\sharp\psi}, \lambda^{\K\!\sH\sharp\psi}) = (\Box(U^{\sS^{op}}\times\mathcal{C}), [\:\Box(U^{\sS^{op}}\times\mathcal{C})]^{H\sharp\psi}, \lambda^{\Box H\sharp\psi}        )
\end{eqnarray*}

\noindent are the corresponding liftings, where $\psi = \Phi(\lambda^{\otimes\widehat{\otimes}})$.

\section{Conclusions and further work}

This article was intended to explore more examples in classic monad theory in order to prove, what the author considers, the \emph{relevance} of the 2-adjunctions of the type $\Adj$-$\Mnd$. This relevance will be significant if the recollection of a numerous quantity of useful examples is done.\\

The development of the article only used the \emph{left} definition, nevertheless, the author hopes that the \emph{right} and the \emph{left-right} case can be completed without any complication whatsoever.\\

As far as further work is concerned, there is a pair of possible connections. The first one, is to take the framework of \emph{multivariable adjunctions} in 
\cite{che_cym} for further analysis using the 2-adjunction and the parametric objects already defined.\\

Second, there might be a further development on categorical duality provided by this parametric objects.\\

\begin{center}
\section*{Acknowledgments}
\end{center}

The author would like to thank to the Consejo Nacional de Ciencia y Tecnolog\'ia (CONACYT) for financial support through the grant SNI-59154. Special thanks, for useful comments on this article to T. Brzezinski, G. Bohm and B. Mesablishvili.\\

The author would like to thank to J. Antonio L. Verver and Fernando Vega for their endless support and patience.\\

$\phantom{br}$

\end{document}